\documentclass[12pt]{article}

\usepackage{amsmath}
\usepackage{t-angles}

\title{Hopf Galois Extension  in Braided Tensor Categories }
\author{
Shouchuan Zhang $^{a,~b}$, \ \   Yao-Zhong Zhang $^b$
  \\ $a$.
  Department  of Mathematics, Hunan University,\\ Changsha
410082, P.R. China \\
$b$. Department of Mathematics, University of Queensland,\\
Brisbane 4072, Australia }
\date{}

\begin{document}
\newtheorem{Theorem}{\quad Theorem}[section]
\newtheorem{Proposition}[Theorem]{\quad Proposition}
\newtheorem{Definition}[Theorem]{\quad Definition}
\newtheorem{Corollary}[Theorem]{\quad Corollary}
\newtheorem{Lemma}[Theorem]{\quad Lemma}
\newtheorem{Example}[Theorem]{\quad Example}
\maketitle \addtocounter{section}{-1}

\begin {abstract} The relation between  crossed product and $H$-Galois
extension in braided tensor category ${\cal C}$ with equivalisers
and coequivalisers  is established.  That is, it is shown that if
there exist an equivaliser and a coequivaliser for any two
morphisms in  ${\cal C}$, then $A = B \# _\sigma H$ is a crossed product
algebra if and only if  the extension
$A/B$ is  Galois, the canonical epic $q: A\otimes A \rightarrow
A\otimes_B A$ is split and $A$ is isomorphic as left $B$-modules
and right $H$-comodules to $B\otimes H$ in ${\cal C}$.

\vskip 0.2cm Keywords: braided Hopf algebra, crossed product
algebra, $H$-Galois extension.
 \end {abstract}

\section {Introduction}
The Hopf Galois extension has its roots in the work of
Chase-Harrison-Rosenberg \cite {CHR65} and  Chase-Sweedler \cite
{CS69}. The general definition about Hopf Galois extension
appeared in \cite {KT81} and the relation between  crossed product
and $H$-Galois extension was obtained in \cite{BM89}\cite
{DT89}\cite{DZ99} for ordinary Hopf algebras. See also the books
by Montgomery \cite {Mo93} and  Dascalescu-Nastasecu-Raianu
\cite{DNR01} for reviews about the main results in this topic.

On the other hand, braided Hopf algebras have attracted much
attention in both mathematics and mathematical physics (see
\cite{AS02}\cite {BD98}\cite{Ke99}\cite{Ma90a}\cite{Ma95b}\cite
{RT93} \cite {Zh03}). In particular, braided Hopf algebras  play an
important role in the classification of finite-dimensional pointed
Hopf algebras (see \cite {AS00}\cite {AS02}). So it is desirable to
generalize the above results to the case of braided tensor
categories. In this paper we show that if there exist an equivaliser
and a coequivaliser for any two morphisms in the braided tensory
category ${\cal C}$, then $A = B \# _\sigma H$ is a crossed product
algebra if and only if  the extension $A/B$ is Galois, the canonical
epic $q: A\otimes A \rightarrow A\otimes_B A$ is split and $A$ is
isomorphic as left $B$-modules and right $H$-comodules to $B\otimes
H$ in ${\cal C}$.

This paper is organized as follows. In section 1  since it is
possible that the category is not concrete, we define the
coinvariants  $B= A^{co H}$ and the tensor product $A\otimes _B A$
over algebra $B$ by equivaliser and coequivaliser for right
$H$-comodule algebra $(A,\psi)$. Our main result in this section
is given in theorem \ref{1.3}. In section 2, we apply the
conclusion in section 1 to  the Yetter-Drinfeld category $^D_D
{\cal YD}$ and give an example to explain our result.

\section { Hopf Galois extension  in braided tensor categories}
In this section we give the relation between  crossed product and
$H$-Galois extension in braided tensor categories.

Let  $({\cal C}, \otimes, I,  C )$ be   a braided tensor category
with equivalisers and coequivalisers, where $I$ is the identity
object and $C$ is the braiding. We denote   $id _W \otimes f$ by
$W \otimes f $ for convenience.  In this section, we assume that
there exist an equivaliser and a coequivaliser for any two
morphisms in  ${\cal C}$.

For two morphisms $f_1 , f_2 : A \rightarrow B,$  set $target (A,
f_1, f_2) =: \{ (D, g) \mid g$  is a morphism from  $D$ to $A$
such that $ f_1 g = f_2g  \}$ and $source  (B, f_1, f_2) =: \{ (D,
g) \mid g$  is a morphism from  $B$ to $D$ such that $ gf_1  =
gf_2  \}$, the   final object of the full subcategory $target (A, f_1,
f_2)$ and the initial object of the full subcategory $source  (B, f_1, f_2)$
are called equivaliser and coequivaliser of $f_1$ and $f_2$ ,
written  as $equivaliser (f_1. f_2)$ and
 $coequivaliser (f_1, f_2)$, respectively .
 Let  $(A, \psi)$ be  a right $H$-comodule algebra and
$(H, m, \eta, \Delta , \epsilon )$  a braided Hopf algebra in
${\cal C}.$ We called the $equivaliser (\psi , id _A \otimes \eta
_H)$ the coinvariants of $A$, written $(A ^{coH}, p).$ We called
the $coequivaliser ( (m_A \otimes A) (A \otimes p \otimes A ) , (A
\otimes m_A) (A \otimes p \otimes A))$ the tensor product of $A$
and $A$ over $A^{co H}$, written $(A \otimes _{A^{co H}} A , q).$
Note that $p $ is monic from $A^{co H}$  to $A$ and $q$ is epic
from $A \otimes A$ to $A\otimes _{A^{co H}} A$ (see \cite {Fa73}).

\begin {Definition}\label {1.1}
Let $(A, \psi )$ be a right $H$-comodule algebra in $  {\cal C}$
and $can '=: (m_A \otimes H)(A \otimes \psi )$ a morphism from $A
\otimes A$ to $A \otimes H$. If there exists an equivalence $can $
in ${\cal C}$ from $A \otimes _{A^{co H}} A$  to $A \otimes H$
such that $can \circ q = can ',$ then we say that  $A$ is a right
$H$-Galois, or the extension $A/A^{co H}$ is Galois.
\end {Definition}

We first recall the crossed product  $B\# _\sigma H$ of $B$ of $H$
in ${\cal C}$  similar to  \cite [Definition 7.1.1]{Mo93}.
$(H,\alpha )$ is said to act weakly on algebra $B $ if the
following conditions are satisfied:

\noindent (WA): \ \  $\alpha (H \otimes m _B) = m_B (\alpha
\otimes \alpha ) (H \otimes C \otimes B) (\Delta _H\otimes B
\otimes B)$ and $\alpha (H \otimes \eta _B) = \eta _B \epsilon _H
. $

\noindent $\sigma $   is called a 2-cocycle from $H\otimes H$  to
$B $  if the following conditions are satisfied:

\noindent (2-COC): \ \  $m _B(\alpha \otimes \sigma ) (H \otimes C
\otimes H) (H \otimes H \otimes \sigma   \otimes m) (H \otimes H
\otimes H\otimes C \otimes H) (\Delta _H \otimes \Delta _H \otimes
\Delta _H)  = m_B (B \otimes \sigma  ) (\sigma  \otimes m \otimes
H )(H \otimes C \otimes H \otimes H) (\Delta _H \otimes \Delta _H
\otimes H) $ and $\sigma (H \otimes \eta _H) = \sigma (\eta _H
\otimes H) = \eta _B \epsilon _H.$

\noindent $(B, \alpha )$   is called a twisted $H$-module   if the
following conditions are satisfied:

\noindent (TM):  $m_B (\alpha \otimes \sigma ) (H \otimes C
\otimes H) (H \otimes H \otimes \alpha \otimes H) (H \otimes H
\otimes H \otimes C) (\Delta _H \otimes \Delta _H \otimes B)= m
_B(B \otimes \alpha ) (\sigma  \otimes m \otimes B)(H \otimes C
\otimes H \otimes B) (\Delta _H \otimes \Delta _H \otimes B)$ and
$\alpha (\eta _H \otimes B)= id _B.$

The crossed product $B \# _\sigma H$ of $B$ and $H$  is  $B
\otimes H$ as an object in ${\cal C},$ with multiplication
\begin {eqnarray*}m_{B\# _\sigma H} &=:  & (m _B \otimes H)( m_B \otimes
\sigma  \otimes m_H)
(B \otimes B \otimes H \otimes C \otimes H)\\
&{ \ }&(B \otimes \alpha \otimes \Delta _H \otimes \Delta _H)(B
\otimes H \otimes C \otimes H)(B \otimes \Delta _H \otimes B
\otimes H)
\end {eqnarray*} and unit $\eta _B \otimes \eta _H.$

\begin {Lemma}\label {1.2}
If $B$  is an algebra and $H$ is a bialgebra in ${\cal C}$, then
$A=B\#_{\sigma} H$   is an algebra with unity element $\eta _A =
\eta _B \otimes \eta _H$
 iff $(H, \alpha )$  acts weakly on $B$ and
$(B, \alpha )$ is a twisted $H$-module with 2-cocycle $\sigma .$
\end {Lemma}

{\bf Proof.} Let  $\hat \sigma =  (\sigma \otimes m_H) (H \otimes
C \otimes H) (\Delta _H \otimes \Delta _H) $ and $\phi _{2,1} =
(\alpha \otimes H) (H \otimes C) (\Delta \otimes B)$. It is clear
that the multiplications in $B\#_\sigma H$ and $B \bowtie _{\phi
_{2,1}}^{\hat \sigma} H $ are the same (see \cite [Proposition
2.2]{BD99}). Consequently, it is sufficient to show that
 \cite [Relation (2.3)]{BD99} holds if and only if $(H, \alpha )$  acts weakly on $B$ and
$(B, \alpha )$ is a twisted $H$-module with 2-cocycle $\sigma .$
 Let us denote the relation (2.3) in  \cite [Proposition 2.2]{BD99} by the following:\\
 (i)
\begin {eqnarray*} \label {e121}
(\alpha \otimes H) (H \otimes C)(\Delta _H \otimes B)(\eta _H
\otimes B)  &=&
(B \otimes \eta _H).   \\
(\sigma \otimes m_H) (H \otimes C \otimes H)(\Delta _H \otimes
\Delta _H)(\eta _H \otimes H)
 &=&
(\eta _B \otimes H).\\
(\alpha \otimes H) (H \otimes C)(\Delta _H \otimes B)(H \otimes
\eta _B)
&=& (\eta _B \otimes H).\\
(\sigma \otimes m_H) (H \otimes C \otimes H)(\Delta _H \otimes
\Delta _H)(H \otimes \eta _H)  &=&
(\eta _B \otimes H). \ \ \ \ \ \ \ \ \ \ \ \ \ \ \ \ \ \ \ \ \ \ \ \   \\
\end {eqnarray*}
(ii)
\begin {eqnarray*}\label {e122}
(\alpha \otimes H) (H \otimes C)(\Delta _H \otimes m_ B) &=&
(m_B\otimes H)(B \otimes \alpha \otimes H) (B \otimes H \otimes C)
\\
&{\ }& (\alpha \otimes \Delta _H \otimes B)
(H \otimes C \otimes B)(\Delta _H \otimes  B\otimes B). \\
\end {eqnarray*}
(iii)
\begin {eqnarray*}\label {e123}
&{ \ }&(m\otimes H)(B\otimes \sigma \otimes m_H) (B \otimes H
\otimes C \otimes H)(\alpha \otimes \Delta _H \otimes \Delta
_H)\\
&{\ }& (H \otimes C \otimes H)   (\Delta_ H \otimes \sigma \otimes
m_H) (H\otimes H \otimes C \otimes H)
 (H \otimes \Delta _H \otimes \Delta _H) \\
  &=& (m_B \otimes
H) (B \otimes \sigma \otimes m_H) (B \otimes H \otimes C \otimes
H)(B \otimes \Delta _H \otimes \Delta _H)\\
&{\ }&(\sigma \otimes m _H \otimes H) (H \otimes C \otimes H
\otimes H)
(\Delta _H \otimes \Delta _H\otimes H).\\
\end {eqnarray*}
(iv)
\begin {eqnarray*}\label {e124}
&{\ }& (m\otimes H)(B\otimes \sigma \otimes m_H) (B \otimes H
\otimes C \otimes H)(\alpha \otimes \Delta _H \otimes \Delta _H)\\
&{\ }& (H \otimes C \otimes H)   (\Delta_ H \otimes \alpha \otimes
H) (H\otimes H \otimes C )(H \otimes \Delta _H \otimes B)\\
 &=&
(m_B \otimes H) (B \otimes \alpha \otimes H) (B \otimes H \otimes
C )(B \otimes \Delta _H \otimes B)(\sigma  \otimes m _H \otimes
B)\\
&{\ }& (H \otimes C \otimes H \otimes B)
(\Delta _H \otimes \Delta _H\otimes B).\\
\end {eqnarray*}
Obviously, if (WA), (2-COC) and (TM) hold then relation (i) holds.
Therefore, we assume that relation (i) holds. Applying $(id _B
\otimes \epsilon _H) $ on relation (iii), we can obtain (2-COC).
It is straightforward that $(ii) \Leftrightarrow $(WA),    $(iii)
\Leftrightarrow $ (2-COC) and $(iv) \Leftrightarrow $(TM).
Consequently, using \cite [Proposition 2.2]{BD99}, we  complete
the proof. \begin{picture}(5,5)
\put(0,0){\line(0,1){5}}\put(5,5){\line(0,-1){5}}
\put(0,0){\line(1,0){5}}\put(5,5){\line(-1,0){5}}
\end{picture}\\


\begin {Theorem}\label {1.3}
Let ${\cal C}$ be a
 braided tensor category  with equivalisers
and coequivalisers for any two morphisms in ${\cal C}$. Then the
following assertions are equivalent:

(i) There exists an invertible 2-cocycle $\sigma : H \otimes H
\rightarrow B, $ and a week action $\alpha $ of $H$ on $B$ such
that $A = B \# _\sigma H$ is a crossed product algebra.

(ii) $(A, \psi )$ is a right $H$-comodule algebra with $B=A ^{co
H}$, the extension $A/B$ is  Galois, the canonical epic $q:
A\otimes A \rightarrow A\otimes_B A$ is split and $A$ is
isomorphic as left $B$-modules and right $H$-comodules to
$B\otimes H$ in ${\cal C}$, where the $H$-comodule operation and
$B$-module operation of $B \otimes H$ are   $\psi _{B \otimes H} =
id _B \otimes  \Delta _H$ and $\alpha _{B \otimes H} = m_B \otimes
id _H$ respectively.
\end {Theorem}
{\bf Proof.} $(i) \Rightarrow  (ii).$ It is clear that $(A, \psi
)$ is a right $B$-comodule algebra under $H$-comodule operation
$\psi = (B\otimes \Delta _H)$. Let  $p =:  id _B \otimes \eta _H$.
For object $D$ and  any morphism $f : D \rightarrow B \# _\sigma
H$ in ${\cal C}$ with $\psi f = (id _A \otimes \eta _H)f$, set
$\bar f =: (id _B \otimes \epsilon ) f$. Thus $p \circ \bar f =
f$, which implies that  $(B, p)$ is the coinvariants of $(A, \psi
).$  Set $\Theta  ' =: (m_A \otimes \eta _B \otimes H) (A \otimes
C_{H, B} \otimes H)(A \otimes H \otimes \sigma ^{-1}\otimes H) (A
\otimes S \otimes S \otimes \Delta )(A \otimes \Delta _H^2) : A
\otimes H \rightarrow A \otimes A$
 and $\Theta =: q \circ \Theta ':
A \otimes H \rightarrow A \otimes _B A$.

 It is remain  to show that $can \circ \Theta  = id$ and $\Theta \circ can = id
.$  For this we first give the relation between the module
operation and the 2-cocycle. Let $\overline{\sigma}$ denote the
convolution-inverse of $\sigma$ and $\sigma^{-1}$ . We denote the
multiplication, comutiplication, antipode , braiding and inverse
braiding by \[
\begin{tangle}
\cu \step[1],\step[2] \cd \step[2], \step[2]\S \step[2],\step[2]\x
\step[2]and\step[2] \xx \step[2], \hbox { respectively .}
\end{tangle}
\]

We first give the relations between the module operation and the
2-cocycle.
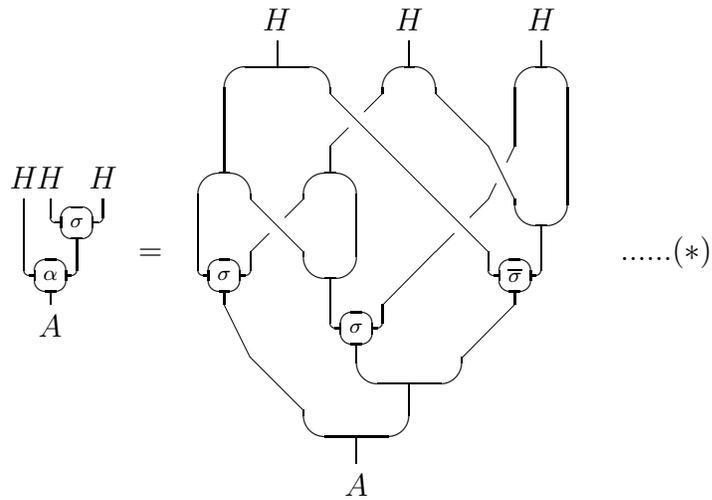
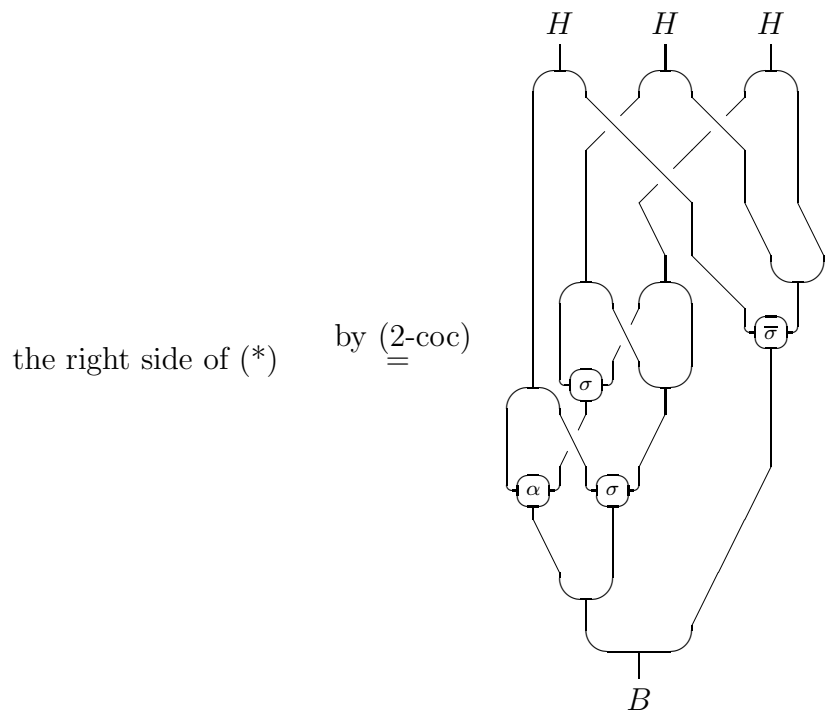
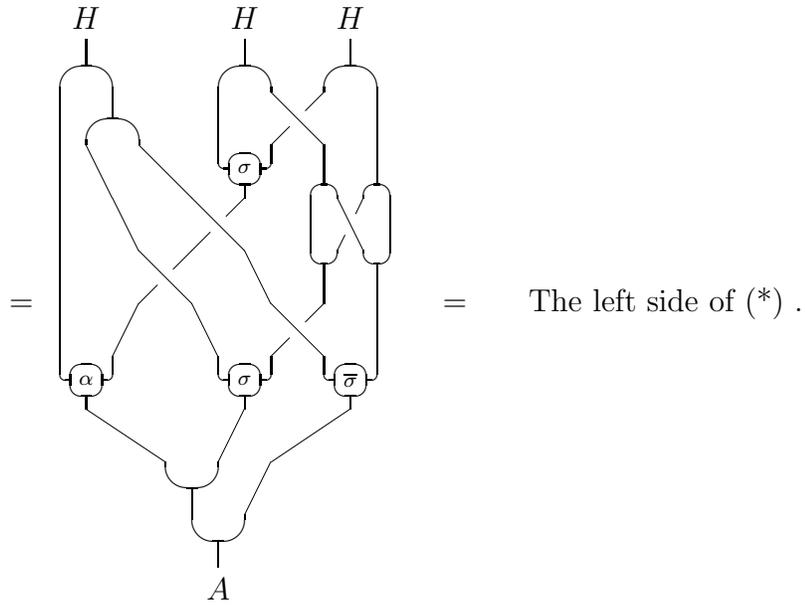
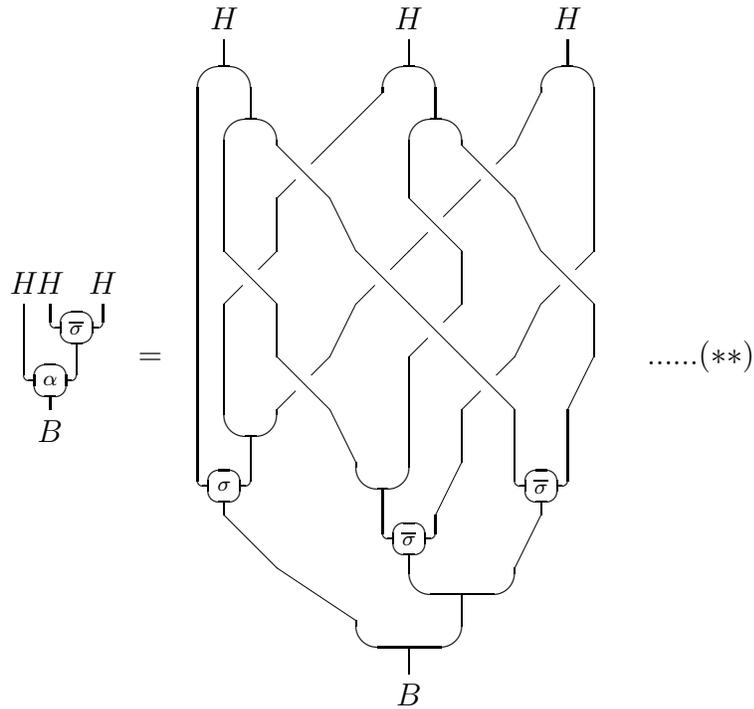
\[
\begin{tangle}
\object{H}\step\object{H}\Step\object{H}\\
\id\step\tu \sigma\\
\tu \alpha\\
\step\object{A}\\
\end{tangle}
\step = \step
\begin{tangle}
\step\Step\object{H}\Step\step\Step\object{H}\Step\Step\step\object{H}\\
\step \Cd\Step\cd\Step\step\cd\\
\step\id\Step\Step\x\Step\nw2\Step\id\Step\id\\
\cd\Step\cd\step\nw2\step\Step\hx\Step\id\\
\id\Step\x\Step\id\Step\step\x\step\cu\\ 
\tu \sigma \Step \cu\Step\ne2\Step\tu {\overline{\sigma}}\\ 
\step\d\Step\step\tu \sigma \Step\Step\ne2\\
\Step\nw2\Step\step\Cu\\
\Step\Step\Cu\\
\Step\Step\Step\object{A}\\
\end{tangle}
\Step  ......(*)
\]
\vskip 0.1cm \noindent
 In fact,
 \vskip 0.1cm
\[ \hbox {
 the right  side  of  (*) } \step \stackrel{ \hbox{ by  (2-coc)}}{=}\step
\begin{tangle}
 \step[2]\object{H} \step[4]\object{H} \step[4]\object{H}\\
 \step \cd  \step[2]\cd \step[2]\cd\\
 \step \id \step[2]\x \step[2]\x \step [2] \id\\
\step \id \step[2] \id \step [2]\x\step [2] \id \step[2] \id\\
\step \id \step[2] \id \step[2]\nw1\step[1]\id \step[2]\nw1\step[1]\nw1\\  
 \step \id\step[1] \cd\step\cd\nw2\step[2]\cu\\  
 \step \id\step[1] \id\Step\hx\Step\id\step[2]\tu
 {\overline{\sigma}}\\  
 \cd\tu \sigma\step \cu\step[3]\id\\  
 \id\step[2]\hx\step[2]\dd\step[2]\step[2]\id\\
\tu \alpha\step\tu \sigma\step[2]\step[2]\dd\\
\step\d\step[2]\id\step[2]\step[2]\dd\\
\step[2]\cu\step[2]\step[1]\dd\\
\step[2]\step[1]\Cu\\
\step[5]\object{B}\\
\end{tangle}
\]
\[
{=}\step
\begin{tangle}
 \step[1]\object{H} \step[6]\object{H} \step[4]\object{H}\\
 \cd \step[2]\step[2]\cd\step[2]\cd\\ 
 \id\step\cd\step[3]\id\step[2]\x\step[2]\id\\ 
 \id\step\d\step\nw2\step[2]\tu \sigma\step[1.5] \hcd\step  \hcd\\  
 \id\step[2]\nw1 \step[2]\x \step[2.5]\id\step\hx\step\id\\ 
 \id\step[3]\x\step[2]\nw1 \step[1.5]\hcu\step\hcu\\  
 \id \step[2]\ne1\step[2]\nw1 \step[2]\x\step[2]\id\\    
 \tu \alpha\step[4]\tu \sigma   \step [2]\tu {\overline{\sigma}}\\
  \step[1]\nw3\step[4] \ne1  \step[3]\ne3\\  
 \step[4]\cu\step \ne1\\  
 \step[5]\cu \\   
 \step[6]\object{A}\\
\end{tangle}
\step[2]
\begin{tangle}
= \hbox {\step[2]The left  side of (*) . }
\end{tangle}
\]
\vskip 0.1cm

\[
\begin{tangle}
\object{H}\step\object{H}\Step\object{H}\\
\id\step\tu {\overline{\sigma}}\\
\tu \alpha\\
\step\object{B}\\
\end{tangle}
\step = \step
\begin{tangle}
 \step[1]\object{H} \step[7]\object{H} \step[6]\object{H}\\
 \cd\step[5]\cd\step[4]\cd\\ 
 \id\step\cd\step[2]\sw2\step[2]\cd\step[2]\ne1\step[2]\id\\  
\id\step\id\step[2]\x\step[3]\id\step[2]\x\step[3]\id\\  
\id\step\id\step[2]\id\step[2]\nw1\step[2]\x\step[2]\nw1\step[2]\id\\  
\id\step\x\step[3]\x\step[2]\id\step[3]\x\\
\id\step\id \step[2]\id \step[2]\ne1 \step[2]\x\step[2]\ne1\step [2]\id \\  
\id\step\id \step[2]\x\step[3]\id \step [2]\x \step[2]\dd\\  
\id\step\cu\step[2]\d\step[2]\id\step[2]\id\step[2]\id\step[2]\id\\  
\tu \sigma \step[4]\cu\step\dd\step[2]\tu {\overline{\sigma}}\\ 
\step\nw2\step[5]\tu {\overline{\sigma}}\step[3]\dd\\  
\step[3]\nw3\step[2]\step[2]\Cu\\ 
\step[6]\Cu\\ 
\step[8]\object{B}\\
 \end{tangle}
 \step[2]......(**)
\]
\vskip 0.1cm \noindent
 Indeed , it is clear that $\alpha(H\otimes
\sigma^{-1})$ is the convolution inverse of $\alpha(H\otimes
\sigma)$  . Therefore , it is sufficient to show that the right side
of $(**)$ is also the convolution inverse of $\alpha(H\otimes
\sigma) .$  See \vskip 0.1cm
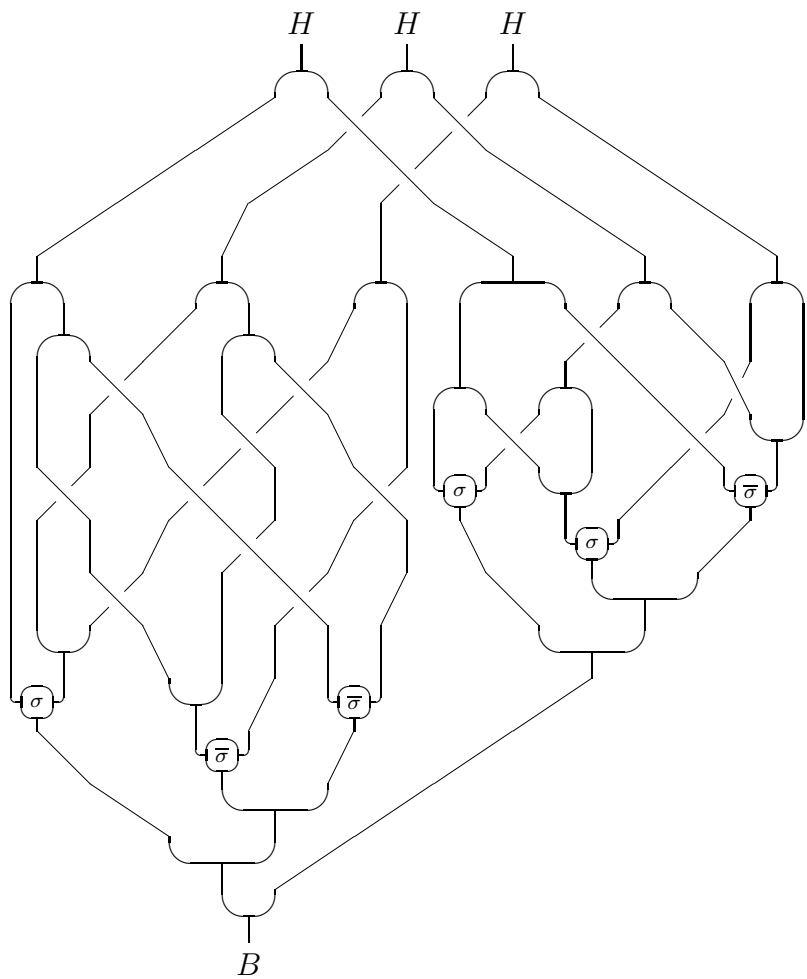
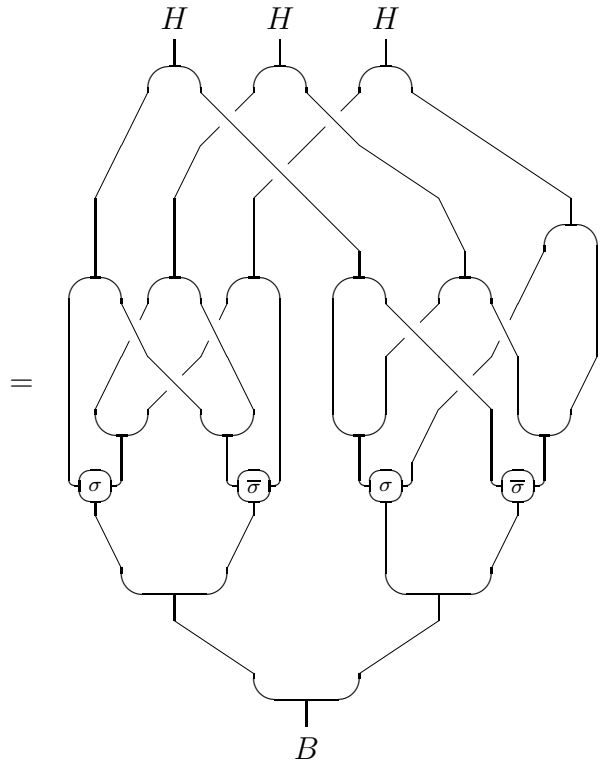
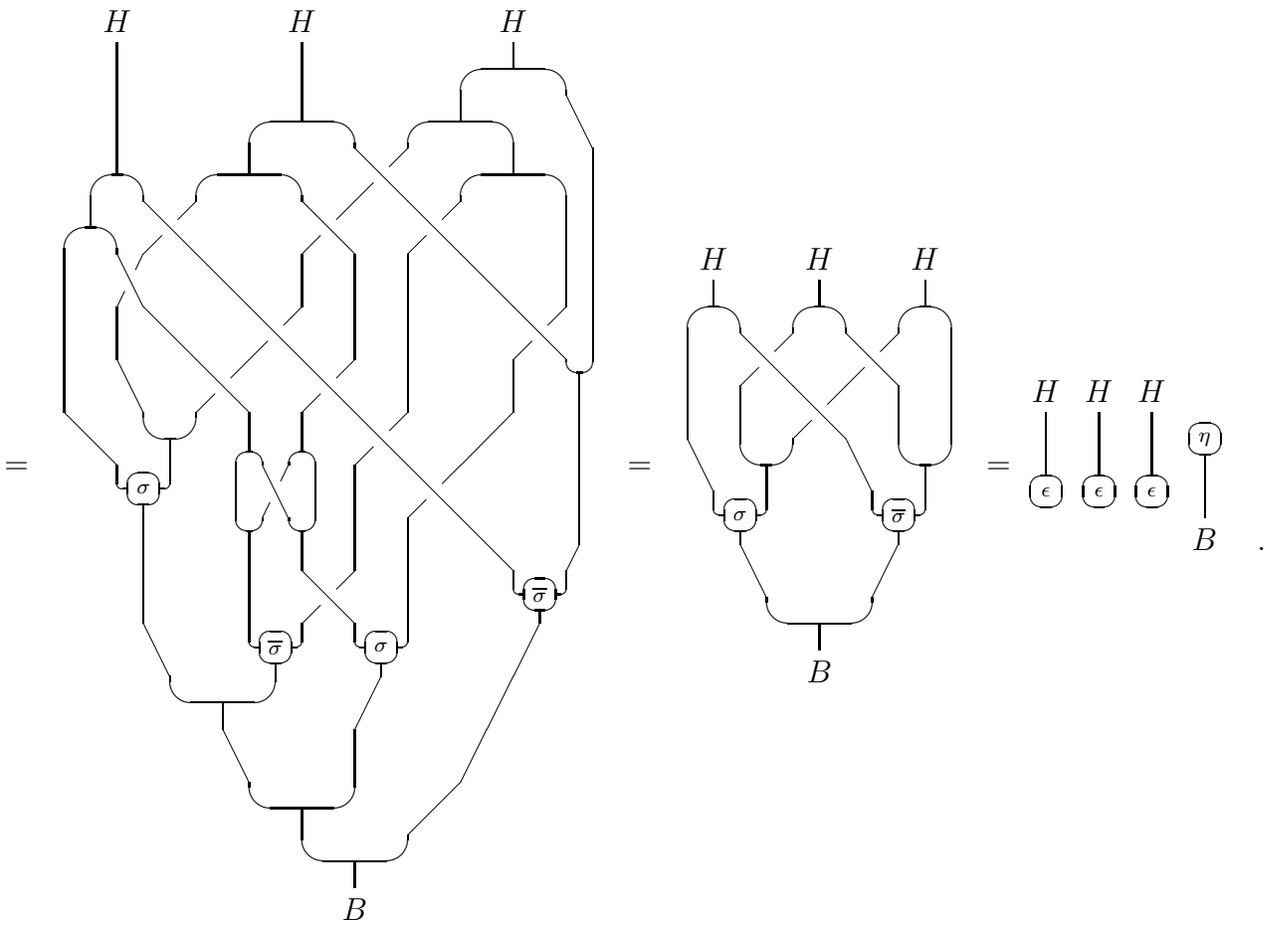
\[
\begin{tangle}
 \step[11]\object{H} \step[4]\object{H} \step[4]\object{H}\\ 
\step[10]\cd\step[2]\cd\step[2]\cd\\ 
\step[9]\ne3\step[2]\x\step[2]\x \step [2]\nw3\\ 
\step[6]\ne3\step[4]\ne3\step[2]\x \step [2]\nw3\step [4]\nw3\\ 
\step[3]\ne3\step[4]\ne1\step[5]\id \step [2]\nw3\step [4]\nw3\step [4]\nw3\\ 
\cd\step[5]\cd\step[4]\cd\step\step \Cd\Step\cd\Step\step\cd\\ 
\id\step\cd\step[2]\sw2\step[2]\cd\step[2]\ne1\step[2]\id \step [2]\id\Step\Step\x\Step\nw2\Step\id\Step\id\\  
\id\step\id\step[2]\x\step[3]\id\step[2]\x\step[3]\id\step [1]\cd\Step\cd\step\nw2\step\Step\hx \Step\id\\  
\id\step\id\step[2]\id\step[2]\nw1\step[2]\x\step[2]\nw1\step[2]\id\step [1]\id\Step\x\Step\id\Step\step\x\step\cu\\  
\id\step\x\step[3]\x\step[2]\id\step[3]\x\step [1]\tu \sigma \Step \cu\Step\ne2\Step\tu {\overline{\sigma}}\\
\id\step\id \step[2]\id \step[2]\ne1 \step[2]\x\step[2]\ne1\step [2]\id \step [2]\d\Step\step\tu \sigma \Step\Step\ne2\\\\  
\id\step\id \step[2]\x\step[3]\id \step [2]\x \step[2]\dd\step [3]\nw2\Step\step\Cu\\  
\id\step\cu\step[2]\d\step[2]\id\step[2]\id\step[2]\id\step[2]\id\step [6]\Cu\\  
\tu \sigma \step[4]\cu\step\dd\step[2]\tu {\overline{\sigma}}\step [7] \ne3\\ 
\step\nw2\step[5]\tu {\overline{\sigma}}\step[3]\dd\step [5] \ne3\\  
\step[3]\nw3\step[2]\step[2]\Cu\step [3] \ne3\\ 
\step[6]\Cu\step [2] \ne3\\ 
\step [8] \cu\\
\step[9]\object{B}\\
 \end{tangle}
\]

\vskip 0.1cm
\[
\step = \step
\begin{tangle}
  \step[4]\object{H} \step[4]\object{H} \step[4]\object{H}\\ 
\step[3]\cd\step[2]\cd\step[2]\cd\\ 
\step[2]\ne1\step[2]\x\step[2]\x \step [2]\nw3\\ 
\step[1]\ne1\step[2]\ne1\step[2]\x \step [2]\nw3\step [4]\nw3\\ 
\step[1]\id\step[3]\id\step[3]\id \step [2]\nw2\step [4]\nw1 \step [3]\cd\\ 
\cd\step\cd\step\cd\step[2]\cd\step[2]\cd\step\dd\step[2]\id\\ 
\id\step[2]\hx\step[2]\hx\step[2]\id \step[2]\id\step[2]
\x\step[2]\hx\step[3]\id\\  
\id\step[1]\ne1 \step[1]\x\step[1]\nw1 \step[1]\id\step[2]
\id \step[2]\id \step [2]\x\step[1]\id\step[2]\ne1\\  
\id\step[1]\cu\step[2]\cu\step[1]\id\step[2]\cu\step \ne1\step[2]
\id\step[1]\cu\\  
\tu \sigma \step[4]\tu {\overline{\sigma}}\step[3] \tu
\sigma \step[2]\step[1]\tu {\overline{\sigma}}\\  
\step[1]\nw1\step[2]\step[2]\ne1\step[2]\step[3]\id\step[2]\step[2]
\dd\\   
\step[2] \Cu\step[6]\Cu\\   
\step[4]\nw3\step[6]\sw3\\
\step[7]\Cu\\
\step[9]\object{B}
 \end{tangle}
\]
\vskip 0.1cm
\[
\step =\step
\begin{tangle}
 \step[2]\object{H} \step[7]\object{H} \step[8]\object{H}\\
\step[2]\id\step[7]\id\step[6]\Cd\\
\step[2]\id\step[5]\Cd\step[2]\Cd\step[2]\d\\  
\step\cd\step[2]\Cd\step[2]\x \step[2]\Cd\step\id\\  
\cd\step[1]\x\step[4]\x\step[2]\x\step[4]\id\step[1]\id\\  
\id\step[2]\hx\step[2]\nw2\step [3]\id\step[2]  \id \step
[2]\id\step[2]\nw2\step[3]\id\step \id\\  
\id\step[2]\id \step[]\nw2\step[3]\x\step[2]\id\step[2]\id\step
[4] \x\step[1]\id\\  
\id\step[2]\nw1\step[2]\x\step[2]\x\step[2]\id\step[4]\id\step
[2] \hcu\\  
\nw2\step[2]\cu\step[1.5]\hcd\step\hcd\step[1.5]\x\step[3]\ne2\step[2.5]\id\\
\step[2]\tu
\sigma\step[2.5]\id\step\hx\step\id\step[1.5]\id \step[2]\x \step[4.5]\id\\  
\step[3]\id\step[3.5]\hcu\step\hcu\step[1.5]\id \step[2]\id\step[2]\nw2 \step [3] \ddh\\   
\step[3]\id\step[2]\step[2]\id\step[2]\x \step[2]\id\step[4]\tu
{\overline{\sigma}} \\ 
\step[3]\d\step[2]\step[1]\tu
{\overline{\sigma}}\step[2]\tu \sigma\step[3]\step\dd\\   
\step[3]\step\Cu\step[3]\dd\step[3]\step[1]\dd\\   
\step[2]\step[3]\step\nw1\step[2]\step[1]\step\id\step[1]\step[3]\dd\\  
\step[7]\Cu\step[2]\sw2\\  
\step[9]\Cu\\  
\step[11]\object{B}
 \end{tangle}
\step =\step
\begin{tangle}
 \step[1]\object{H} \step[4]\object{H} \step[4]\object{H}\\
 \cd\step[2]\cd\step[2]\cd\\
 \id\step[2]\x\step[2]\x\step[2]\id\\
\id\step[2]\id \step[2]\x\step[2]\id\step[2]\id\\
 \nw1\step[1]\cu\step[2]\nw1\step[1]\cu\\
 \step \tu \sigma\step[2]\step[2]\tu {\overline{\sigma}}\\
 \step[2]\d\step[4]\dd\\
 \step[3]\Cu\\
 \step[5]\object{B}\\
 \end{tangle}
 \step =\step
\begin{tangle}
\object{H} \step[2]\object{H}\step[2]\object{H}\\
\id\step[2]\id\step[2]\id  \step[2]\Q \eta\\
 \QQ \epsilon\step[2]\QQ \epsilon  \step[2]\QQ \epsilon     \step[2]\id\\
 \step[2]\step[4]\object{B} \step [2].
 \end{tangle}
\step[6]
\]


\vskip 0.1cm Now we see that
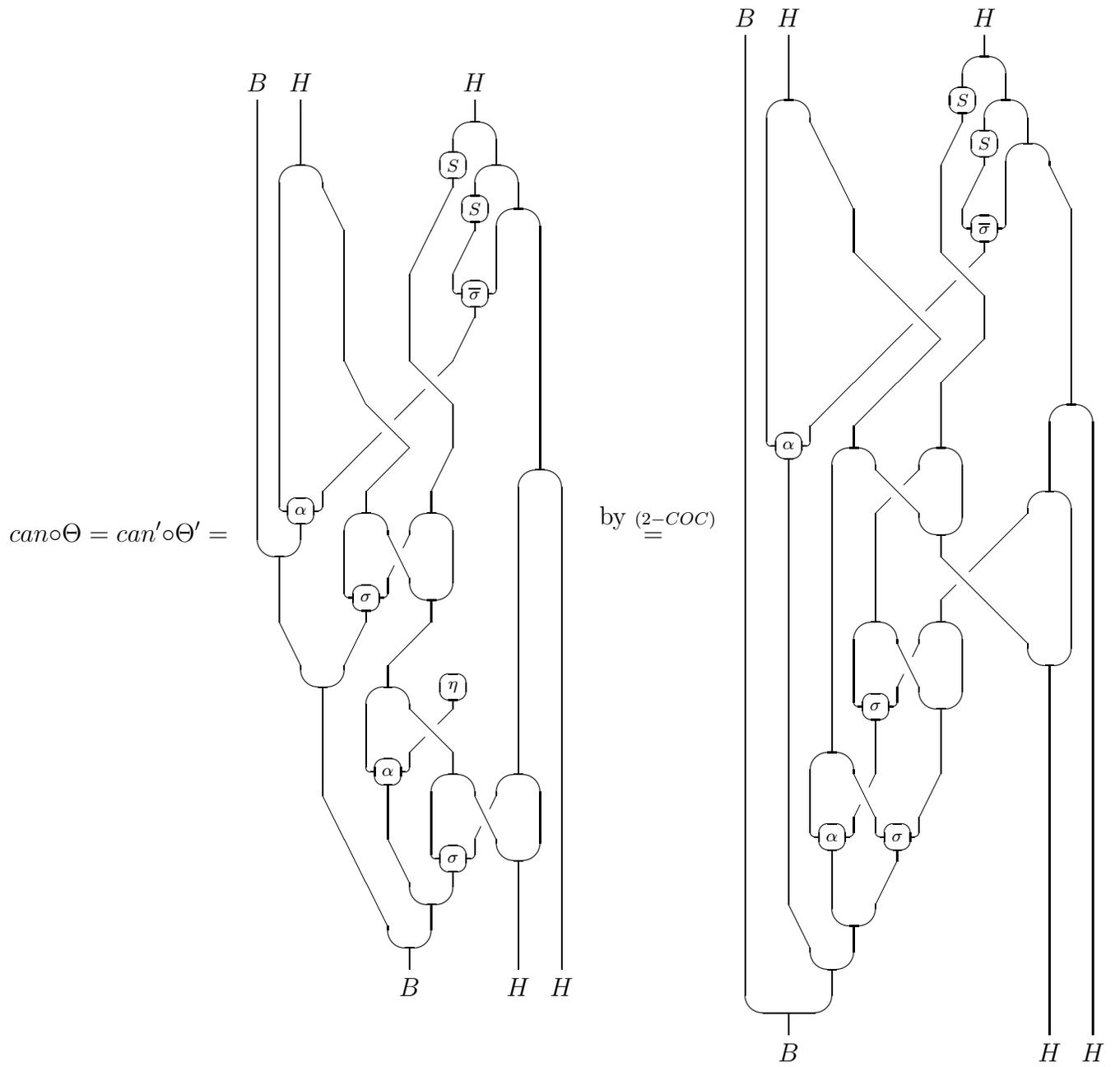
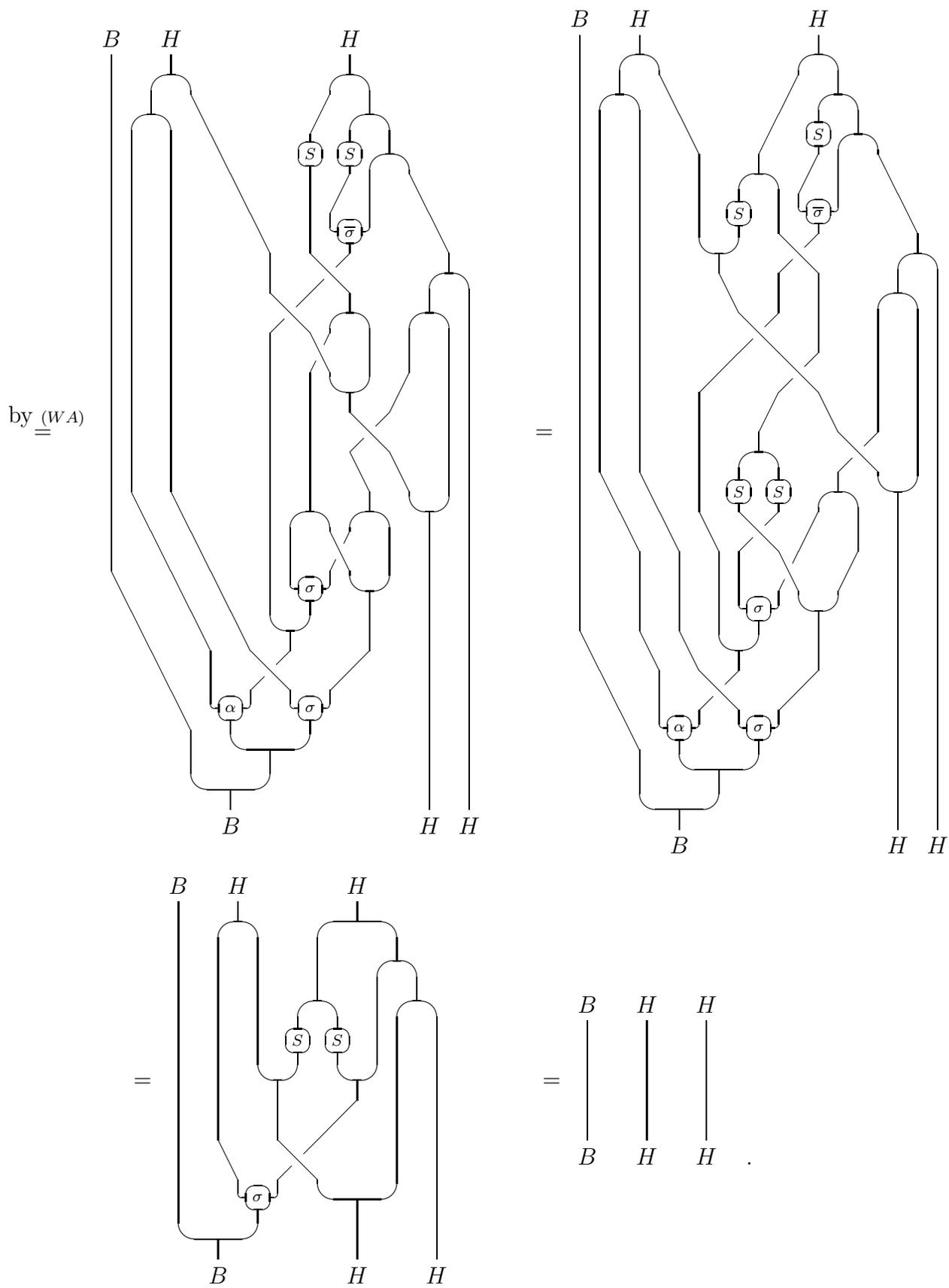
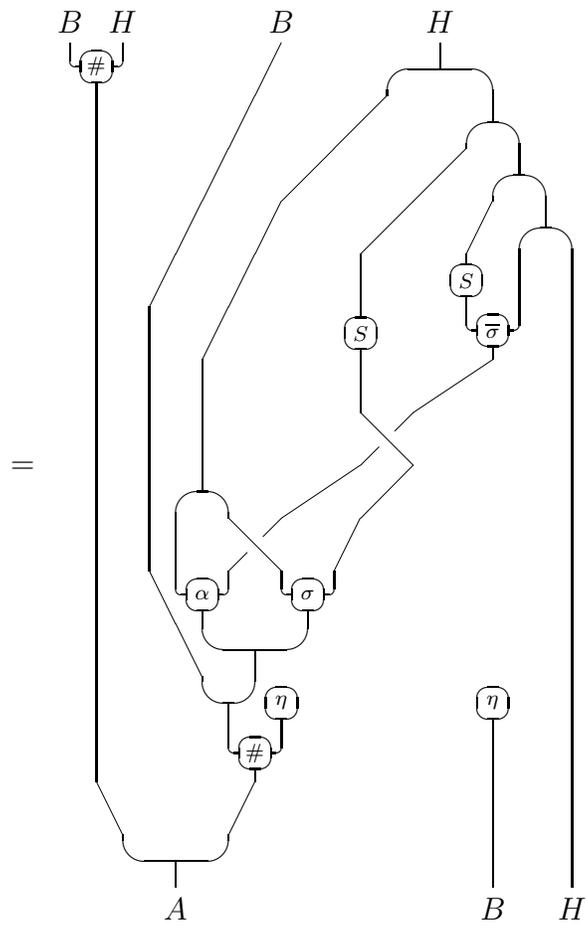
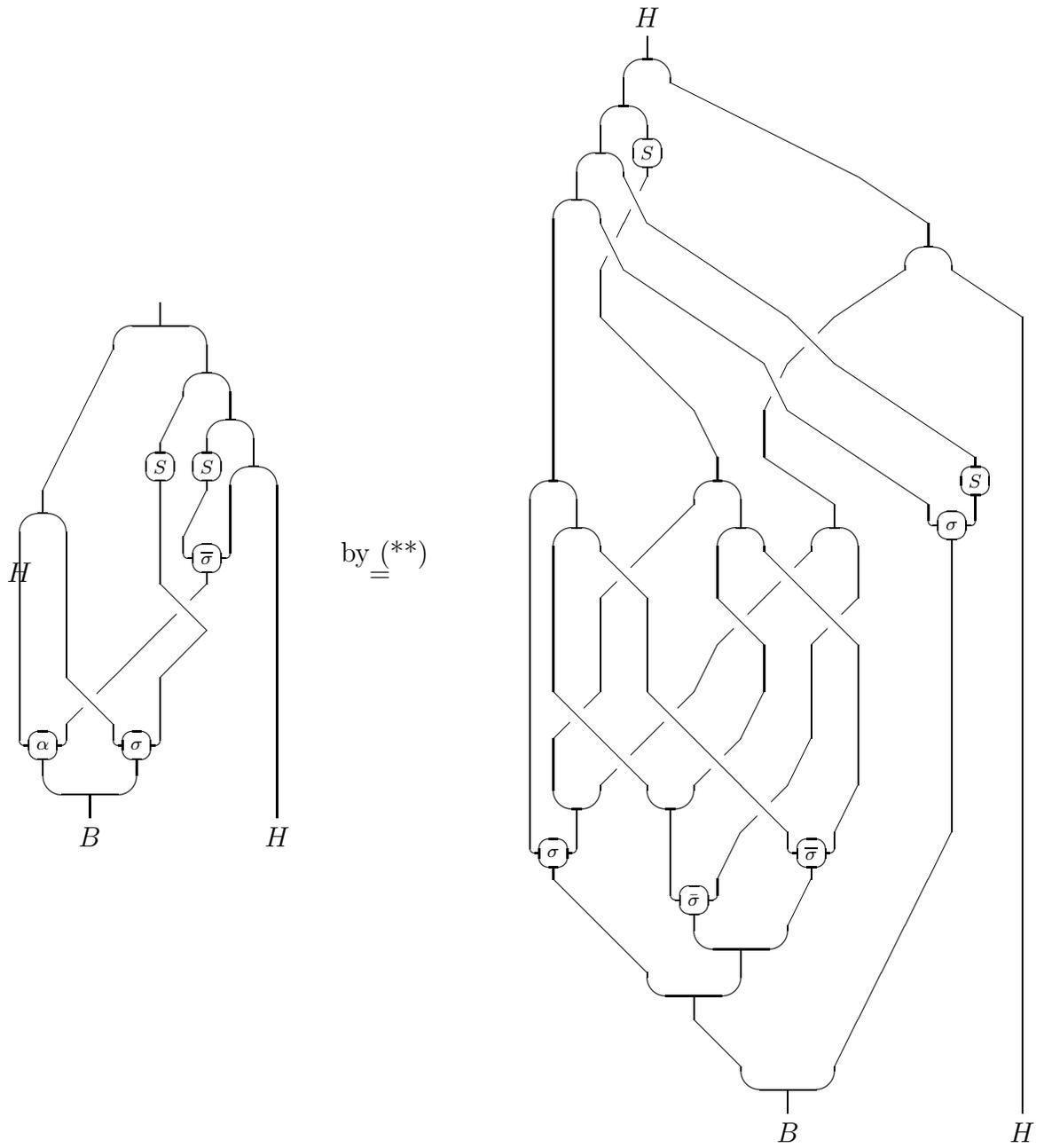
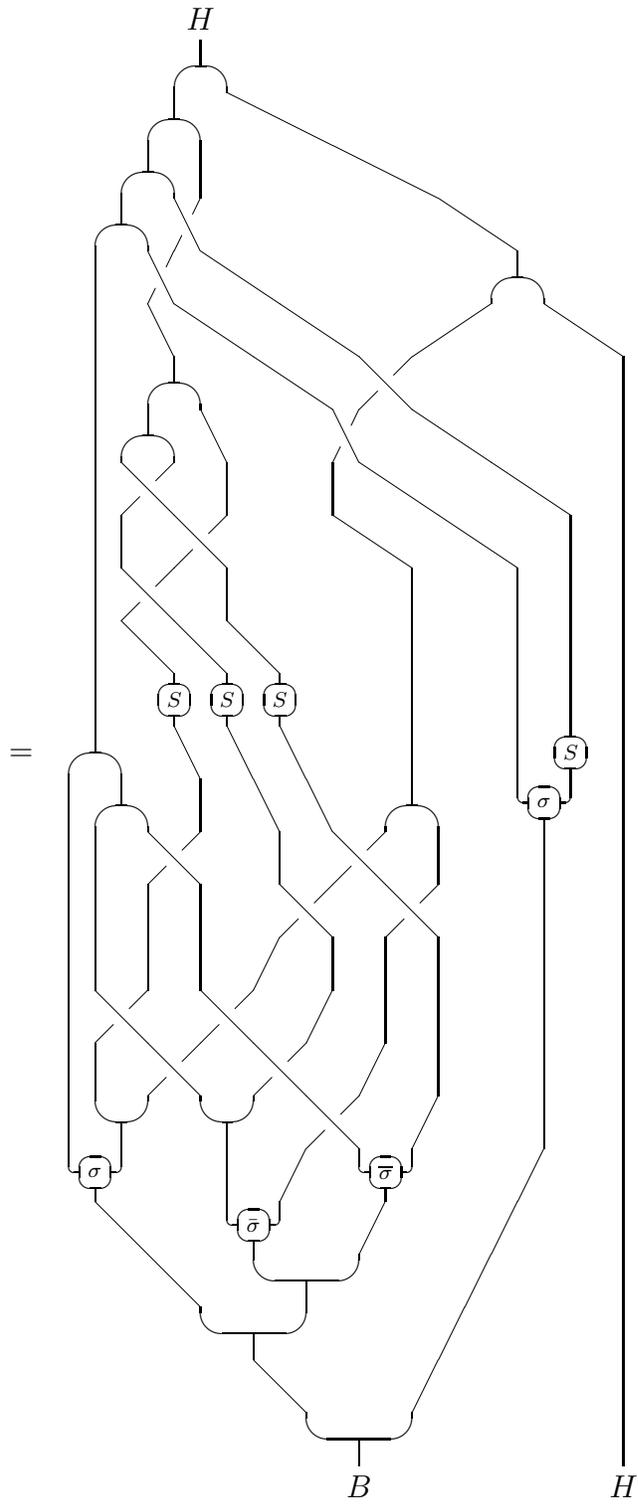
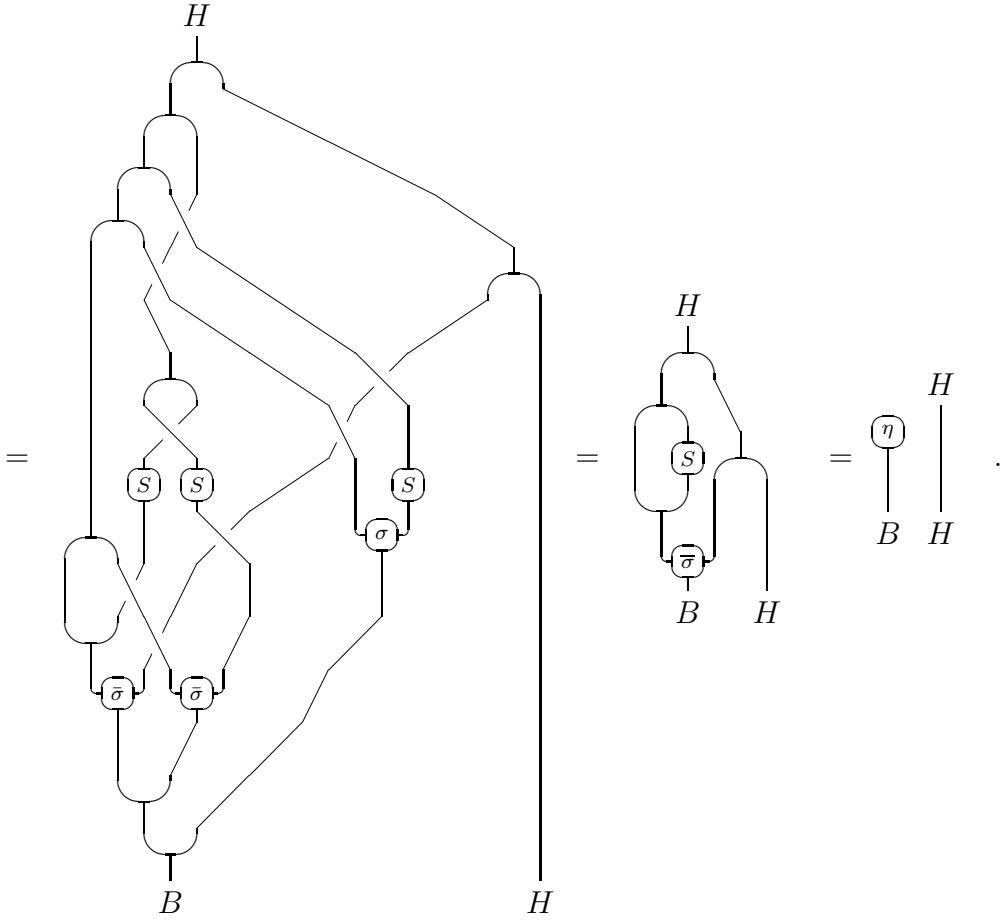
\[
 can \circ \Theta  = can '\circ \Theta  '=\step
\begin{tangle}
\object{B}\step[2]\object{H}\step[8]\object{H}\\
\id\step[2]\id\step[7]\cd\\  
\id\step\cd\step[2]\step[2]\step[2]\S\step\cd\\ 
\id\step\id\step[2]\d\step[2]\step[2]\dd\step\S\step\cd\\  
\id\step\id\step[2]\step\id\step[2]\step\dd\step[1]\dd\step\id\step[2]
\id\\   
\id\step\id\step[2]\step\id\step[3]\id\step[2]\tu
{\overline{\sigma}}\step[2]\id \\ 
\id\step\id\step[2]\step\id\step[2]\step[1] \id\step[2] \ne1\step[3]\id\\
\id\step\id\step[2]\step\nw1\step[2] \x\step[4]\id\\
\id\step\id\step[4]\x\step [2]\id\step[4]\id\\

\id \step \id \step[3]\ne2\step[1]\ne2\step [1]\ne1\step[3]\cd\\

\id\step\tu \alpha\step \cd\step\cd\step[]\step[2]\id\step[2]\id\\

\cu\step[2]\id\step[2]\hx\step[2]\id\step[2]\step
\id\step[2]\id\\  

\step \id \step [3]\tu \sigma  \step  \cu\step\step[2]\id\step[2]\id\\  

\step \nw1   \step [2]\ne1 \step[2]\ne2\step[2] \step[2]\id\step[2]\id\\   

\step\step\cu\step\cd\step[2]\Q \eta \step[2]\step\id\step[2]\id\\   

\step\step[2]\id\step[2] \id\step[2]\x\step[2]\step\id\step[2]\id\\    
\step[2]\step\id\step[2]\tu
\alpha\step\cd\step\cd\step\id\\  
\step[2]\step\d\step[2]\id\step[2]\id\step[2]\hx\step[2]\id\step\id\\ 
\step[2]\step[2]\d\step\d\step\tu \sigma\step\cu\step\id\\
\step[2]\step[2]\step\d\step\cu\step[2]\step\ \id\step[2]\id\\
\step[2]\step[2]\step[2]\cu\step[2]\step[2]\id\step[2]\id\\
\step[7]\object{B}\step[5]\object{H}\step[2]\object{H}\\
 \end{tangle}
 \step  \stackrel { \hbox { by } (2-COC) } {=}\step
\begin{tangle}
\object{B}\step[2]\object{H}\step[9]\object{H}\\
\id\step[2]\id\step[7]\step\cd\\  
\id\step\cd\step[2]\step[2]\step[2]\step\S\step\cd\\ 
\id\step\id\step[2]\d\step[2]\step[2]\step\dd\step\S\step\cd\\  
\id\step\id\step[2]\step\d\step[2]\step\step\id\step[1]\dd\step\id\step[2]
\d\\   
\id\step\id\step[2]\step[2]\id\step[2]\step[2]\id\step\tu
{\overline{\sigma}}\step[2]\step\id \\ 
\id\step\id\step[2]\step\step\nw2\step[1]\step[2]\x \step[3]\step\id\\
\id\step\id\step[6]\x\step[2]\id\step[2]\step[2]\id\\
\id\step\id\step[5]\ne2\step[1]\ne2\step[1] \ne2 \step[4]\id\\
\id\step\id\step[3]\ne2\step[1]\ne2\step[2] \id \step[5]\cd\\
\id\step\tu
\alpha\step\cd\step[2]\cd\step[2]\step[2]\id\step[2]\id\\  
\id\step[2]\id\step[2]\id\step[2]\x\step[2]\id\step[2]\step\cd\step\id\\  
\id\step[2]\id\step[2]\id\step[2]\id\step[2]\cu\step\sw2\step[2]\step\id\step\id\\ 
\id\step[2]\id\step[2]\id\step[2]\id\step[2]\step
\x\step[2]\step[2]\id\step\id\\  
\id\step[2]\id\step[2]\id\step\cd\step\cd\step\nw2\step[2]\step\id\step\id\\  
\id\step[2]\id\step[2]\id\step\id\step[2]\hx\step[2]\id\step[2]\step\cu\step\id\\
\id\step[2]\id\step[2]\id\step\tu
\sigma\step\cu\step[2]\step[2]\id\step[2]\id\\ 
\id\step[2]\id\step\cd\step\id\step[2]\step\id\step[2]\step[2]\step\id\step[2]\id
\\  
\id\step[2]\id\step\id\step[2]\hx\step[2]\dd\step[3]\step[2]\id\step[2]\id\\   
\id\step[2]\id\step\tu \alpha\step\tu
\sigma\step[2]\step[2]\step[2]\id\step[2]\id\\    
\id\step[2]\id\step[2]\id\step[2]\dd\step[2]\step[2]\step[2]\step\id\step[2]\id\\
\id\step[2]\d\step\cu\step[8]\id\step[2]\id\\  
\id\step[2]\step\cu\step[9]\id\step[2]\id\\ 
\Cu\step[10]\id\step[2]\id\\  
\step[2]\object{B}\step[12]\object{H}\step[2]\object{H}\\
 \end{tangle}
 \step\step  \step
\]

\vskip 0.1cm
\[
\step \stackrel { \hbox { by } (WA) }  {=} \step
\begin{tangle}
\object{B}\step[3]\object{H}\step[9]\object{H}\\
\id\step[2]\cd\step[7]\cd\\ 
\id\step\cd\step\d\step[2]\step[2]\step\dd\step\cd\\   
\id\step\id\step[2]\id\step[2]\d\step[2]\step[2]\S\step[2]\S\step\cd\\ 
\id\step\id\step[2]\id\step[2]\step\d\step[2]\step\id\step\dd\step\id\step[2]\d\\  
\id\step\id\step[2]\id\step[2]\step[2]\d\step[2]\id\step\tu
{\overline{\sigma}}\step[2]\step\d\\  
\id\step\id\step[2]\id\step[2]\step[2]\step[]\id\step[2]\x\step[2]\step[2]\cd\\  
\id\step\id\step[2]\id\step[2]\step[2]\step\x\step\cd\step[2]\cd\step\id\\  
\id\step\id\step[2]\id\step[2]\step[2]\step\id
\step[2]\hx\step[2]\id\step[2]\id\step[2]\id\step\id\\  
\id\step\id\step[2]\id\step[2]\step[2]\step\id\step[2]\id\step\cu\step[1]\ne1\step[2]\id
\step\id\\     
\id\step\id\step[2]\id\step[2]\step[2]\step\id\step[2]\id\step[2]\x\step[2]\step
\id\step\id \\   
\id\step\id\step[2]\id\step[2]\step[2]\step\id\step[2]\id\step[2]\nw1
\step [1]\nw1\step[]\step
\id\step\id \\   
\id\step\d\step\d\step[2]\step[2]\id\step\cd\step\cd\step\cu\step\id\\  
\id\step[2]\d\step\d\step[2]\step\id\step\id\step[2]
\hx\step[2]\id\step[2]\id\step[2]\id\\  
\d\step[2]\d\step\d\step[2]\id\step \tu \sigma\step\cu\step[2]
\id\step[2]\id\\ 
\step\d\step[2]\d\step\d\step\cu\step[2]\step\id\step[2]\step
\id\step[2]\id\\  
\step[2]\d\step[2]\id \step[2]\x\step [3]
\ne2\step[1]\step[2]\id\step[1]\step\id\\   
\step\step[2]\d\step\tu \alpha\step[2]\tu \sigma
\step[4]\step\id\step[2]\id\\   
\step[2]\step[2]\id\step[]\step\Cu\step[2]\step[4]\id\step[2]\id\\ 
\step[2]\step[2]\Cu\step[2]\step[2]\step[2]\step[2]\id\step[2]\id\\  
\step[6]\object{B}\step[10]\object{H}\step[2]\object{H}\\
 \end{tangle}
 \step\step
\step=\step
\begin{tangle}
\object{B}\step[3]\object{H}\step[9]\object{H}\\
\id\step[2]\cd\step[7]\cd\\ 
\id\step\cd\step\d\step[2]\step[2]\step\dd\step\cd\\   
\id\step\id\step[2]\id\step[2]\d\step[2]\step\dd\step[2]\S\step\cd\\ 
\id\step\id\step[2]\id\step[2]\step\id\step[2]\cd\step\dd\step\id\step[2]\d\\  
\id\step\id\step[2]\id\step[2]\step\id\step[2]\S\step[2]\id\step\tu
{\overline{\sigma}}\step[2]\step\d\\  
\id\step\id\step[2]\id\step[2]\step\cu\step[2]\x \step[2]\step[2]\cd\\  
\id\step\id\step[2]\id\step[3]\step\nw1\step[2]\id  \step[2]\id \step[3] \cd \step [1]\id \\  

\id\step\id\step[2]\id\step[2]\step[3]\x\step[2]\id\step\step[2]\id
\step[2] \id
\step\id\\  
\id\step\id\step[2]\id\step[4]\ne2 \step[2]\x\step[2]\step\id\step[2]\id\step\id\\  
\id\step\id\step[2]\id\step[2]\step\id\step[3]\ne1\step[2]\nw1\step[2]\id\step[2]\id
\step\id\\     
\id\step\id\step[2]\id\step[2]\step\id\step[2]\cd\step[3]\x
\step[2]
\id\step\id\\   
\id\step\d\step\d\step[2]\id\step[2]\S\step[2]\S\step[2]\cd\step\cu\step\id\\  
\id\step[2]\d\step\d\step\d\step\x\step\sw1\step[2]\id\step[2]\id\step[2]\id\\  
\id\step[3]\id\step[2]\id\step[2]\id\step[]\id
\step[2]\hx\step[2]\sw1\step[2]\id\step[2]
\id\\  
\id\step[3]\id\step[2]\id\step[2]\id\step[] \tu
\sigma\step\cu\step[3]
\id\step[2]\id\\ 

\d\step[2]\d\step\d\step\cu\step[2]\step\id\step[2]\step[2]
\id\step[2]\id\\  

\step[]\d\step[2]\id \step[2]\x
\step[2]\step[]\ne2\step[2]\step[2]\id\step[2]\id\\   

\step[2]\d\step\tu \alpha\step[2]\tu \sigma
\step[5]\step\id\step[2]\id\\   

\step[2]\step[]\id\step[2]\Cu\step[2]\step[5]\id\step[2]\id\\ 

\step[1]\step[2]\Cu\step[2]\step[2]\step[2]\step[3]\id\step[2]\id\\  
\step[5]\object{B}\step[11]\object{H}\step[2]\object{H}\\
 \end{tangle}
  \step[2]
\]
\vskip 0.1cm
\[
\step=\step
\begin{tangle}
\object{B}\step[3]\object{H}\step[6]\object{H}\\
\id\step[2]\cd\step[2]\step\Cd\\
\id\step[2]\id\step[2]\id\step[2]\step\id\step[2]\step\cd\\
\id\step[2]\id\step[2]\id\step[2]\cd\step[2]\id\step\cd\\
\id\step[2]\id\step[2]\id\step[2]\S\step[2]\S\step[2]\id\step\id\step[2]\id\\
\id\step[2]\id\step[2]\cu\step[2]\cu\step\id\step[2]\id \\
\id\step[2]\id\step[2]\step\id \step[2]\step[1]\ne2\step[2]\id\step[2]\id\\
\id\step[2]\nw1\step[2]\x \step[2]\step[2]\id\step[2]\id\\
\id\step[3]\tu \sigma \step[]\step\Cu\step[2]\id\\
\Cu\step[2]\step[2]\step[]\id\step[2]\step[2]\id\\
\step[2]\object{B}\step[7]\object{H}\step[4]\object{H}\\
 \end{tangle}
 \step[2]\step[2]
\step=\step
\begin{tangle}
\object{B}\step[3]\object{H}\step[3]\object{H}\\
\id\step[3]\id\step[3]\id\\
\id\step[3]\id\step[3]\id\\
\id\step[3]\id\step[3]\id\\
\object{B}\step[3]\object{H}\step[3]\object{H} \step [2].\\
 \end{tangle}
 \]
\vskip 0.1cm
\[
\step\Theta '\circ  can ' \step =\step
\begin{tangle}
\object{B}\step[2]\object{H}\step[2]\object{B}\step[3]\object{H}\\
\tu \# \step[2]\id\step[2]\cd\\
\step\id\step[2]\step\tu \# \step\cd\\
\step\Cu\step[2]\S\step\cd\\
\step[2]\step\id\step[2]\step\dd\step\S\step\cd\\
\step[2]\step\id\step[2]\step\id\step\dd\step\id\step[2]\id \\
\step[2]\step\id\step[2]\step\id\step\tu
{\overline{\sigma}}\step[2]\id\\
\step[2]\step\id\step[2]\step\x\step[2]\step\id\\
\step[2]\step\id\step[2]\step\tu \#\step\Q \eta\step\step\id\\
\step\step[2]\Cu\step[2]\id\step\step\id\\
\step[5]\object{A}\step[4]\object{B}\step[2]\object{H}\\
 \end{tangle}
\step =\step
\begin{tangle}
\object{B}\step[2]\object{H}\step[6]\object{B}\step[6]\object{H}\\
\tu \# \step[2]\step[2]\step\dd\step[2]\step[2]\Cd\\  
\step\id\step[2]\step[2]\step\dd\step[2]\step\sw2\step[2]\step[2]\cd\\  
\step\id\step[2]\step[2]\dd\step[2]\sw2\step[2]\step[2]\sw2\step[2]\cd\\  
\step\id\step[2]\step\dd\step[2]\dd\step[2]\step[2]\cd\step\dd\step\cd\\  
\step\id\step[2]\dd\step[2]\dd\step[2]\step[2]\step\id\step[2]
\id\step\S\step[2]\id\step[2]\id\\   
\step\id\step[2]\id\step[2]\cd\step[2]\step[2]\step\x\step\tu
{\overline{\sigma}}\step[2]\id \\  

\step\id\step[2]\id\step[2]\id\step[2]\d\step[2]\step[2]\id \step
[2]
\x\step[3]\id\\   

\step\id\step[2]\id\step[2]\id\step[3]\nw2\step[]\step[2]
\x\step[2]\id\step[3]\id\\   

\step\id\step[2]\id\step[2]\id\step[3]\step[2]
\x\step[2]\S\step[2]\S
\step[2]\step\id\\ \\   

\step\id\step[2]\id\step[2]\id\step[2]\step[2] \ne2\step[2] \id
\step [2] \id
\step[2]\id \step[2]\step\id\\ \\   

\step\id\step[2]\id\step\cd\step[2]\id\step[3]\ne2
\step[1]\ne2\step[1]\ne2\step[3]\id\\   
\step\id\step[2]\id\step\id\step[2]\x\step[2]\x\step[2]\id
\step[2]\step[2]\step\id\\    
\step\id\step[2]\d\tu \alpha\step[2] \tu \sigma \step[2]\cu
\step[2]\step[2]\step\id\\   
\step\id\step[2]\step\d\Cu\step[2]\sw2\step\step[2]\step[2]\step[2]\id\\  
\step\id\step[2]\step[2]\cu\step\sw3\step[2]\step[2]\step[2] \step
\Q \eta\step[2]\step\id\\  
\step\id\step[2]\step[2]\step\tu \#
\step[2]\step[2]\step[2]\step[2]
\id\step[2]\step\id\\   
\step\d\step[2]\step[2]\dd\step[2]\step[2]\step[2]\step[2]\step
\id
\step[2]\step\id\\   
\step[2]\Cu\step[10]\id\step[2]\step\id\\  
\step[4]\object{A}\step[12]\object{B}\step[3]\object{H}\\
 \end{tangle}
 \step[4]
\]

\[
\step =\step
\begin{tangle}
\object{B}\step[2]\object{H}\step[6]\object{B}\step[6]\object{H}\\
\tu \#\step[2]\step[2]\step\dd\step[2]\step[2]\Cd\\  
\step\id\step[2]\step[2]\step\dd\step[2]\step\sw2\step[2]\step[2]\cd\\  
\step\id\step[2]\step[2]\dd\step[2]\sw2\step[2]\step[2]\sw2\step[2]\cd\\  
\step\id\step[2]\step\dd\step[2]\dd\step[2]\step\sw2\step\step\step\dd\step\cd\\  
\step\id\step[2]\dd\step[2]\dd\step[2]\step[2]\id\step\step[2]
\step\S\step[2]\id\step[2]\id\\   
\step\id\step[2]\id\step[2]\dd\step\step[2]\step[2]\S\step\step[2]\step\tu
{\overline{\sigma}}\step[2]\id\\  
\step\id\step[2]\id\step[2]\id\step[2]\step[2]\step[2]\id\step[2]
\sw3\step\step[2]\step[2]\id\\   
\step\id\step[2]\id\step[2]\id\step[2]\step[4] \x \step[6]\id\\   
\step\id\step[2]\id\step\cd\step[4]\ne3\step[]\ne2
\step[2]\step[3]\step\id \\   
\step\id\step[2]\id\step\id\step[2]\x\step[2]\dd\step\step[2]
\step[2]\step[2]\step\id\\    
\step\id\step[2]\d\tu \alpha\step[2] \tu \sigma \step[2]\step[2]
\step[2]\step[2]\step\id\\   
\step\id\step[2]\step\d\Cu\step[2]\step\step\step[2]\step[2]\step[2]\id\\  
\step\id\step[2]\step[2]\cu\step\Q
\eta\step\step[2]\step[2]\step[2] \step
\Q \eta\step[2]\step\id\\  
\step\id\step[2]\step[2]\step\tu
\#\step[2]\step[2]\step[2]\step[2]
\id\step[2]\step\id\\   
\step\d\step[2]\step[2]\dd\step[2]\step[2]\step[2]\step[2]\step
\id
\step[2]\step\id\\   
\step[2]\Cu\step[10]\id\step[2]\step\id\\  
\step[4]\object{A}\step[12]\object{B}\step[3]\object{H}\\
 \end{tangle}
\]
\vskip 0.1cm \noindent and
\[\step[6]\object{H}
\begin{tangle}
\step[2]\step[2]\Cd\\
\step[2]\step\dd\step[2]\step\cd\\
\step[2]\dd\step[2]\step \dd\step\cd\\
\step\dd\step[2]\step[2]\S\step[2]\S\step\cd\\
\cd\step[2]\step[2]\id\step\dd\step\id\step[2]\id\\
\id\step[2] \id\step[2]\step[2]\id\step\tu
{\overline{\sigma}}\step[2]\id\\
\id\step[2]\id\step[4]\x\step[2]\step\id\\
\id\step[2]\id\step[3]\ne2\step[]\ne2\step[3]\id\\
\id\step[2]\x\step[2]\id\step[2]\step[2]\step\id\\
\tu \alpha\step[2]\tu \sigma\step[2]\step[2]\step\id\\
\step\Cu\step[2]\step[2]\step[2]\id\\
\step[3]\object{B}\step[8]\object{H}\\
 \end{tangle}
 \step\step \stackrel {\hbox { by (**)} }{=} \step[4]
\begin{tangle}
 \step[5]\object{H} \\
\step[2]\step\step\cd\\  
\step[2]\step\cd\step\nw4\\  
\step[2]\cd\step\S\step[2]\step[2]\step\nw4\\   
\step\cd\step\hx\step[9]\nw3\\   
\step\id\step[2]\hx\step[1]\nw3\step[10]\cd\\  
\step\id\step[2]\id\step[1]\nw3\step[3]\nw3\step [6]\ne3\step[2]\nw3\\  
\step\id\step[2]\nw2\step[3]\nw3\step [3]\x\step[8]\id\\  
\step\id\step[4]\nw2\step[4]\hx\step[2]\nw3\step[7]\id\\  
\step\id\step[6]\nw1\step[2]\id\step[1]\nw3\step[4]\nw3\step[4]\id\\  
 \cd\step[5]\cd\step[]\nw3\step[3]\nw3\step[4]\S\step[2]\id\\ 
 \id\step\cd\step[2]\sw2\step[2]\cd\step[2]\cd\step[2]\step\tu \sigma\step[2]\id\\  
\id\step\id\step[2]\x\step[3]\id\step[2]\x\step[2] \id\step[2]
\step[2]\id\step[3]\id\\  
\id\step\id\step[2]\id\step[2]\id\step[3]\x\step[2] \x\step[2]
\step[2]\id\step[3]\id \\  
\id\step\id\step[2]\id\step[2]\id\step[2]\ne1\step[2] \id\step[2]
\id\step[2]\id \step[4]\id \step[3]\id \\  
\id\step\x\step[2]\x\step[2]\ne1\step[2]\id\step[2]\id\step[2]\step[2]
\id\step[3]\id\\  
\id\step\id \step[2]\x\step[2]\x\step[2]\ne1\step[2]\id \step[3]
\step\id\step[3]\id\\  
\id\step\cu\step[2]\cu\step[2]\x\step[2]\ne1
\step[4]\id\step[3]\id
\step[2]\step[2]\step\\ 
\tu \sigma \step[4]\id\step[2] \ne1\step[2]\tu {\overline{\sigma}}
\step[2]\step[2]\dd\step[3]\id\\ 
\step\nw2\step[4]\tu {\bar {\sigma}}
\step[3]\ne1\step[4]\dd\step[3]\step\id\\  
\step[3]\nw2\step[2]\step[1]\Cu\step[2]\step[2]\dd\step[2]\step[3]\id\\  
\step[5]\Cu\step[2]\step[2]\step\dd\step[2]\step[3]\step\id\\  
 \step[7]\nw2\step[2]\step[2]\step\dd\step[2]\step[2]\step[3]\id\\ 
\step[9]\Cu\step[8]\id\\  
\step[11]\object{B}\step[10]\object{H}\\\end{tangle}
\]

\vskip 0.1cm

\[ \step = \step
\begin{tangle}
 \step[5]\object{H} \\
\step[2]\step\step\cd\\  
\step[2]\step\cd\step\nw4\\  
\step[2]\cd\step[]\id \step[2]\step[2]\step\nw4\\   
\step\cd\step\hx\step[9]\nw3\\   
\step\id\step[2]\hx\step[1]\nw3\step[10]\cd\\  
\step\id\step[2]\nw1\nw3\step[3]\nw3\step [6]\ne3\step[2]\nw3\\  
\step\id\step[2]\cd\step[2]\nw3\step [3]\x\step[8]\id\\  
\step\id\step[1]\cd\step[1]\nw1\step[4]\hx\step[2]\nw3\step[7]\id\\  
\step\id\step[1]\x\step[2]\id \step[4]\id\step[1]\nw3\step[4]\nw3\step[4]\id\\  
\step\id\step[1]\id \step [2]\x\step[2]\step[2]\nw3\step[3]\nw3\step[4]\id\step [2]\id\\  
\step\id\step[1]\x \step [2]\id \step [2]\step[2] \step[3]\id\step[4]\id\step[2]\id\step[2]\id\\  
\step\id\step[1]\nw2 \step [1]\nw2 \step[1]\nw2 \step [2]\step[1] \step[3]\id\step[4]\id\step[2]\id\step[2]\id\\  
\step\id\step[3]\S \step [2]\S\step[2]\S \step[1]\step[1] \step[3]\id\step[4]\id\step[2]\id\step[2]\id\\  
\cd\step[2]\nw1\step[]\nw1\step[1]\nw1\step[4]\id \step [4] \id \step [2]\S\step[2]\id\\ 
 \id\step\cd\step[2]\id \step[2]\nw1\step[1]\nw1 \step [2]\cd\step[2]\step\tu \sigma\step[2]\id\\  
\id\step\id\step[2]\x\step[3]\id\step[2]\x\step[2] \id\step[2]\step[2]\id\step[3]\id\\  
\id\step\id\step[2]\id\step[2]\id\step[3]\x\step[2] \x\step[2]\step[2]\id\step[3]\id \\  
\id\step\id\step[2]\id\step[2]\id\step[2]\ne1\step[2] \id\step[2]
\id\step[2]\id \step[4]\id \step[3]\id \\  
\id\step\x\step[2]\x\step[2]\ne1\step[2]\id\step[2]\id\step[2]\step[2]
\id\step[3]\id\\  
\id\step\id \step[2]\x\step[2]\x\step[2]\ne1\step[2]\id \step[3]
\step\id\step[3]\id\\  
\id\step\cu\step[2]\cu\step[2]\x\step[2]\ne1
\step[4]\id\step[3]\id
\step[2]\step[2]\step\\ 
\tu \sigma \step[4]\id\step[2] \ne1\step[2]\tu {\overline{\sigma}}
\step[2]\step[2]\dd\step[3]\id\\ 
\step\nw2\step[4]\tu {\bar {\sigma}}
\step[3]\ne1\step[4]\dd\step[3]\step\id\\  
\step[3]\nw2\step[2]\step[1]\Cu\step[2]\step[2]\dd\step[2]\step[3]\id\\  
\step[5]\Cu\step[2]\step[2]\step\dd\step[2]\step[3]\step\id\\  
 \step[7]\nw2\step[2]\step[2]\step\dd\step[2]\step[2]\step[3]\id\\ 
\step[9]\Cu\step[8]\id\\  
\step[11]\object{B}\step[10]\object{H}\\\end{tangle}
\]

\vskip 0.1cm

\[ \step = \step
\begin{tangle}
 \step[5]\object{H} \\
\step[2]\step\step\cd\\  
\step[2]\step\cd\step\nw4\\  
\step[2]\cd\step[]\id \step[2]\step[2]\step\nw4\\   
\step\cd\step\hx\step[9]\nw3\\   
\step\id\step[2]\hx\step[1]\nw3\step[10]\cd\\  
\step\id\step[2]\nw1\nw3\step[3]\nw3\step [6]\ne3\step[2]\id \\  
\step\id\step[2]\cd\step[2]\nw3\step [3]\x\step[5]\id\\  
\step\id\step[2]\x\step[2] \step[3]\hx\step[2]\id\step[4]\step[1]\id\\  
\step\id\step[2]\S \step [2]\S\step[2]\step[2]\ne3\step[1]\id \step[2]\S \step [5]\id\\  
\cd\step[1]\id \step[2]\x \step [4]\tu \sigma \step[2] \step[3]\id\\  
\id\step[2]\hx\step [1]\ne1 \step[2]\id  \step [2] \step[3]\id \step[6]\id\\  
\cu\step[1]\hx\step [2]\ne1 \step[4]\ne2 \step [2]\step[1] \step[3]\id\\  
\step [1]\tu {\bar \sigma } \step[]\tu {\bar \sigma } \step[3]
\ne1\step[8] \id\\
\step[2]\id\step[2]\ne1\step[3]\ne2\step[8]\step\id\\
\step[2]\cu\step[2]\ne2\step[2]\step[8]\step\id\\  
\step[3]\cu\step[2]\step[2]\step[3]\step[2]\step[3]\step\id\\  
\step[4]\object{B}\step[14]\object{H}\\\end{tangle} \step = \step
\begin{tangle}
\step[2]\object{H}\\
\step\cd\\
\cd\step\d\\
\id\step[2]\S\step\cd\\
\cu\step\id\step[2]\id\\
\step\tu {\overline{\sigma}}\step[2]\id\\
\step[2]\object{B}\step[3]\object{H}
 \end{tangle}
 \step\step = \step
\begin{tangle}
\step[2]\object{H}\\
\Q \eta  \step[2]\id\\
\id\step[2]\id\\
 \object{B}\step[2]\object{H}\\
 \end{tangle}
  \step\step .\step[2]
\]
\vskip 0.1cm \noindent
 Consequently ,
\[
\Theta ' \circ can '\step=\step
\begin{tangle}
\object{B}\step[2]\object{H}\step[2]\object{B}\step[6]\object{H}\\
\id\step[2]\id\step[2]\id\step[2]\Q \eta\step[2]\Q
\eta\step[2]\id\\
\tu \#\step[2]\tu \#\step[2]\id\step[2]\id\\
\step\Cu\step[2]\step\id\step[2]\id\\
\step[2]\step\object{A}\Step\Step\step\object{B}\Step\object{H}\\
 \end{tangle}
 \]
\vskip 0.1cm \noindent See that
\begin {eqnarray*}
\Theta\circ can \circ q \ &=& \ q \circ \Theta '\circ can '\\
&=& q(m_{A}\otimes A)(A\otimes p\otimes A)(A\otimes B\otimes
\eta_{B} \otimes H)\\
&=& q(A\otimes m_{A})(A\otimes p\otimes A)(A\otimes B\otimes
\eta_{B} \otimes H)\\
&=&  q  \ .
\end {eqnarray*}
\noindent Thus $\Theta \circ can = id _{A \otimes_{B} A}$.
Furthermore, since $\Theta = q \circ \Theta '$ we have $q \circ
(\Theta ' \circ \Theta ^{-1})= id $, which implies $q$ is a split
epic morphism.

 (ii) $\Rightarrow $  (i). Since $q$ is split, there exists morphism $w : A\otimes _B A \rightarrow
 A \otimes A$ such that
 $q\circ w = id _{A\otimes _B A}$. Let $\Theta ' =: w \circ can ^{-1}: \ A\otimes H
  \rightarrow A \otimes  A$, $can ' =: can \circ q :\  A\otimes A \rightarrow A\otimes H$
 and
 $\gamma =: \Phi (\eta _A \otimes H): H \rightarrow
A,$
 $u = can ^{-1}(\eta _A \otimes H):  H \rightarrow A \otimes _B A$,
  $u' = \Theta ' (\eta _A \otimes H):
 H \rightarrow A \otimes A$ and
$\mu = m _A(A \otimes p \otimes \epsilon _H) (A \otimes \Phi
^{-1})u ' : H \rightarrow A$. Obviously, $can ^{-1} = q \circ
\Theta '$, \ \    $can '\circ  u' = can \circ q \circ \Theta '
(\eta _A \otimes H ) = can \circ can ^{-1} (\eta _A \otimes H) =
(\eta _A \otimes H)$ and $\psi \circ \gamma = (\gamma \otimes H)
\Delta $.
\begin {eqnarray}\label {e1}(A \otimes
\Theta  ) u' = (u' \otimes H)\Delta \end {eqnarray} $\hbox { since
} (can ' \otimes H)(A \otimes \psi )u' = (m_A \otimes H \otimes
H)(A \otimes \psi \otimes H) (A \otimes \psi  ) u' = (m_A \otimes
\Delta ) (A \otimes  \psi )u'  = \eta _A \otimes \Delta
   \hbox { and }  (can ' \otimes  H) (u' \otimes H) \Delta =
(can '\circ  u' \otimes H)\Delta = \eta _A \otimes \Delta .$ { }
\begin {eqnarray}\label {e2}  m_A u' = \eta _A \epsilon _H
\end {eqnarray}
 $ \hbox {since } m u ' = (m_A \otimes
\epsilon _H) (A \otimes \psi  ) u' = (A \otimes \epsilon ) can'
\circ  u' = (A \otimes \epsilon _H)(\eta _A \otimes H) =\eta _A
\epsilon _H. $ {}
\begin {eqnarray}\label {e3} (m_A
\otimes A )(A \otimes u')\psi = \eta _A \otimes A \end {eqnarray}
$ \hbox { since }  can ' (m_A \otimes A) (A \otimes u') \psi =
(m_A \otimes H) (m_A \otimes \psi) (A \otimes u') \psi = (m_A
\otimes H) (m_A \otimes A \otimes H) (A \otimes u' \otimes H)(A
\otimes \Delta _H ) \psi  = \psi \hbox { \ \ by (\ref {e1}), (\ref
{e2}),   \ \ and \ \  }  can ' (\eta _A \otimes A ) = \psi.$


Now we show that $\mu $ is the convolution  inverse of $\gamma$.
Indeed,
\begin  {eqnarray*}  \gamma * \mu &=& m  (\Phi  \otimes \mu ) (\eta _B \otimes \Delta _H)\\
&=&m _A (A \otimes \mu )\psi \Phi (\eta _B \otimes H) \hbox {\ \
since } \Phi  \hbox { is an } H
\hbox {-comodule morphism }\\
&=&(m_A \otimes \epsilon _H)(A \otimes m_A \otimes H)(A \otimes A
\otimes p \otimes H) (A \otimes A \otimes \Phi ^{-1})
(A \otimes u')\psi  \Phi (\eta _B \otimes H)\\
&=&(p \otimes \epsilon _H)(\eta _B \otimes H)  \ \ \hbox  { by  (\ref {e3})}\\
&=&\eta _A \otimes \epsilon _H \hbox {\ \ \ \ \ \ \ \ \ \ and }\\
 \mu * \gamma &=&  m _A (m_A \otimes A ) (A \otimes p  \otimes \epsilon \otimes  A )
 (A \otimes \Phi ^{-1} \otimes \Phi  ) (u' \otimes \eta _A \otimes H)\Delta _H \\
&=& m _A (A \otimes \Phi ) (A \otimes B   \otimes \epsilon
\otimes H) (A \otimes \Phi ^{-1}
\otimes H) (u' \otimes H) \Delta _H  \\
 &{ \ \ }&  \hbox {\ \
since } \Phi  \hbox { is a } B
\hbox {-module morphism }\\
&= & m _A (A \otimes \Phi ) (A \otimes B   \otimes \epsilon
\otimes H) (A \otimes \Phi ^{-1}
\otimes H) (A \otimes \psi )u'  \hbox {\ \ by (\ref {e1})}\\
&=& m _A (A \otimes \Phi ) (A \otimes B   \otimes \epsilon
\otimes H) (A \otimes B \otimes \Delta _H)
(A \otimes \Phi ^{-1})u' \\
&{ \ \ }&  \hbox {\ \ since } \Phi ^{-1} \hbox { is an } B
\hbox {-comodule morphism }\\
&=& m_A (A \otimes \Phi) (A \otimes \Phi ^{-1}) u'\\
&=& m_A u'\\
&=& \epsilon _H \eta _A \hbox {\ \ by (\ref {e2})} .
 \end {eqnarray*}
 Thus $\gamma $ has a convolution inverse $\mu.$ Since both
$\psi \gamma ^{-1}$ and $C_{H, A} (S \otimes \gamma ^{-1})\Delta $
are the convolution inverse of $\psi \gamma,$ $\psi \gamma ^{-1} =
C_{H, A} (S \otimes \gamma ^{-1})\Delta .$  Set
\begin  {eqnarray*} \alpha ' &=:& m_A  (A  \otimes m_A ) (\gamma \otimes p \otimes \gamma ^{-1}) (H \otimes
C_{H, B}) (\Delta \otimes B):  \ \ H \otimes B \rightarrow A, \\
\sigma '&=:& m_A  (m_A \otimes \gamma ^{-1} ) (\gamma  \otimes
\gamma \otimes m)(H \otimes C \otimes H)
 (\Delta _H \otimes \Delta _H): \ \ H \otimes H \rightarrow A.\\
\omega '&=:& m_A  (\gamma \otimes m_A) (m_H \otimes C_{A, A})(H
\otimes H
 \otimes \gamma ^{-1} \otimes \gamma ^{-1})
(H \otimes C \otimes H) \\
&{\ }& (\Delta _H \otimes \Delta _H): \ \ H \otimes H \rightarrow A.\\
 \end {eqnarray*}
Since $ \psi m _A (p \otimes p )=  (A \otimes \eta _H ) m _A (p
\otimes p ),$  there exists $m_B: B \otimes B \rightarrow B$ such
that $p m_B = m_A (p \otimes p).$ Similarly, there exist $\eta _B
: I \rightarrow B,$  $\sigma : H \otimes H \rightarrow B,$ $\alpha
: H \otimes B \rightarrow B $ and $\omega : H \otimes H
\rightarrow B$ such that $p \eta _B = \eta _A$, $p \sigma  =
\sigma '$,  $p \alpha = \alpha '$ and $p \omega =\omega '$. It is
easy to check that $\omega $ is the convolution inverse of $\sigma
$. Furthermore, $(B, m_B, \eta _B)$ is an algebra in ${\cal C}.$

 Now we show that conditions (2-$COC$), $(WA)$ and (TM) hold.

\vskip 0.1cm
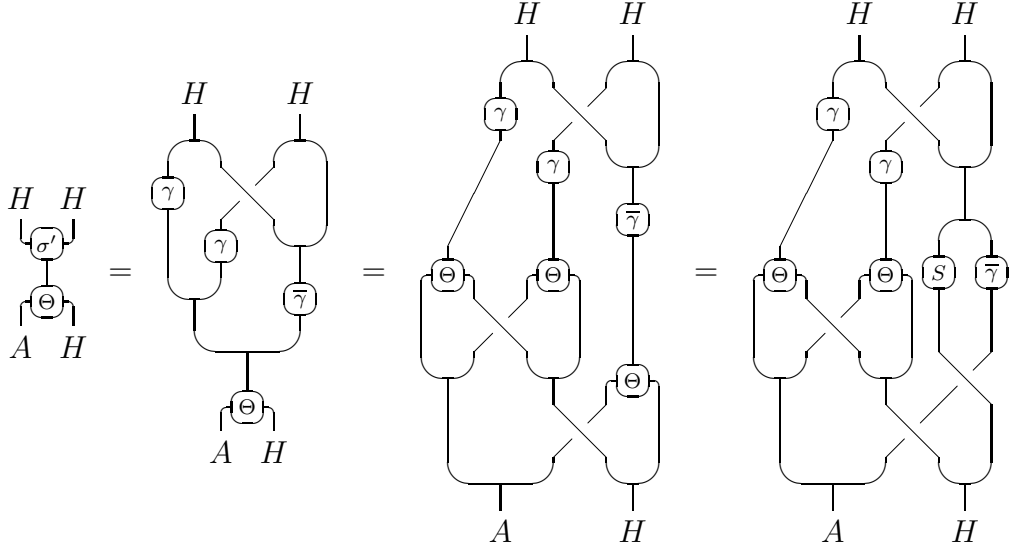
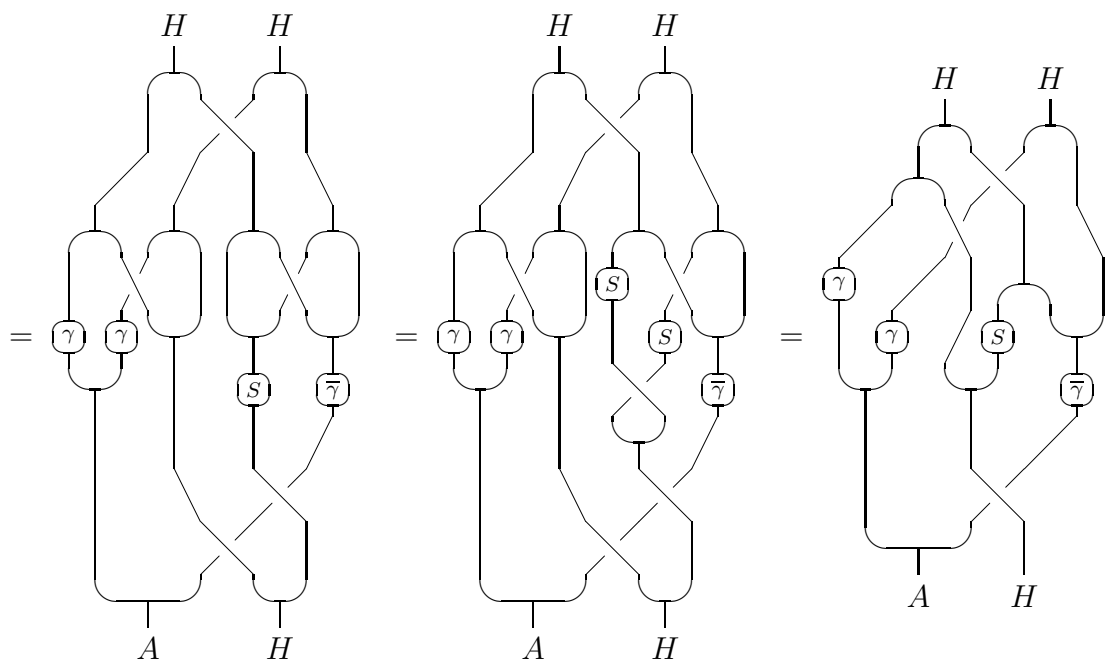
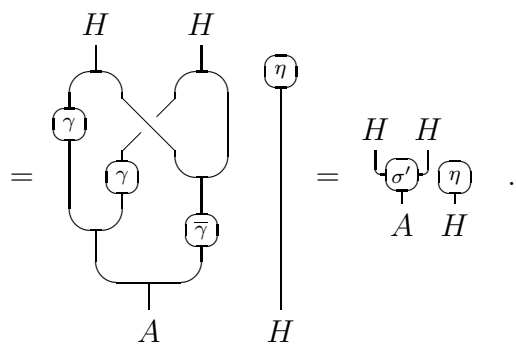
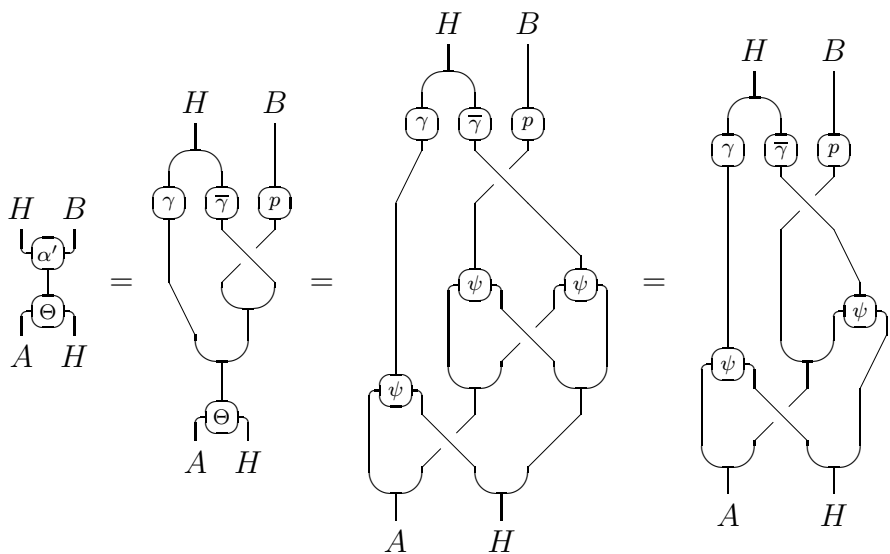
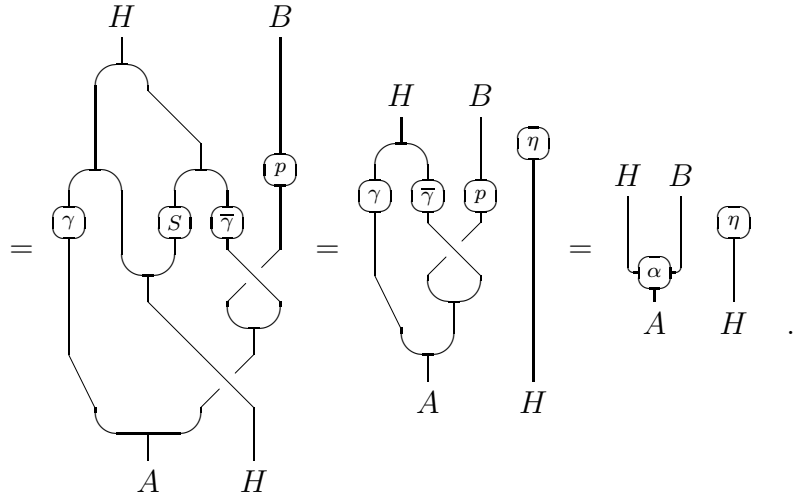
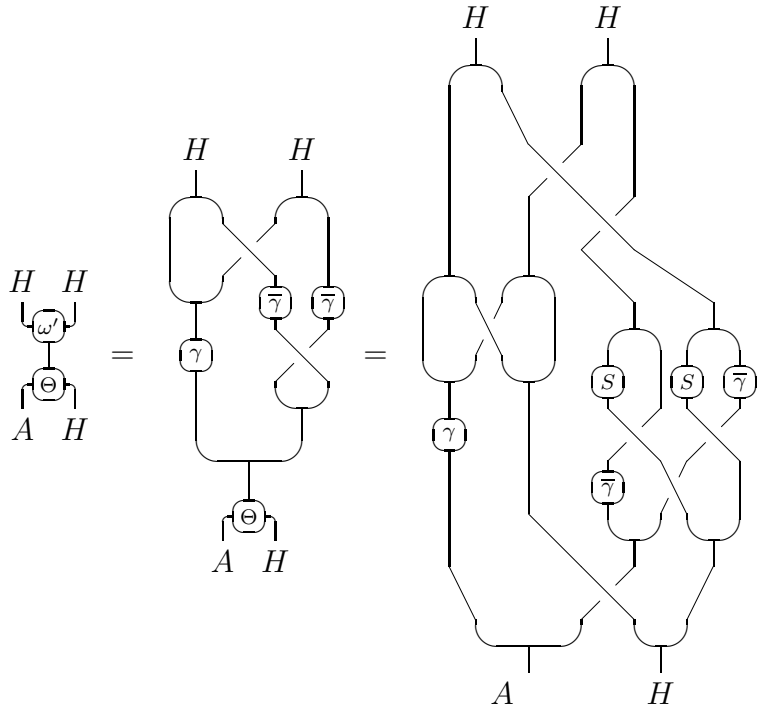
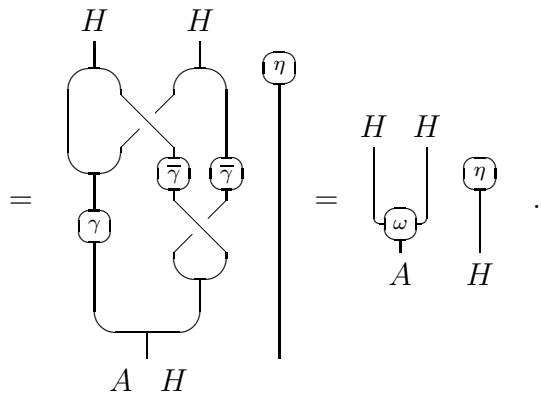
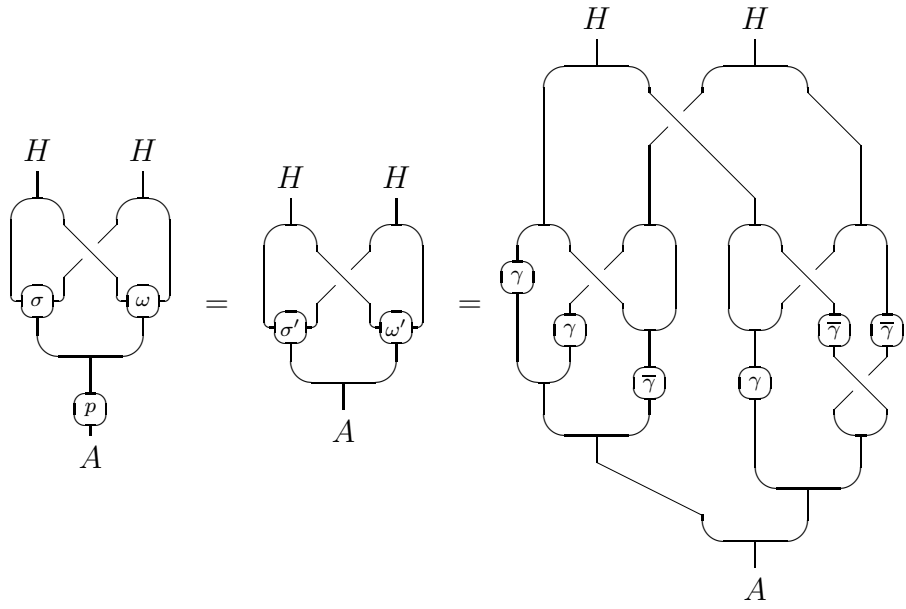
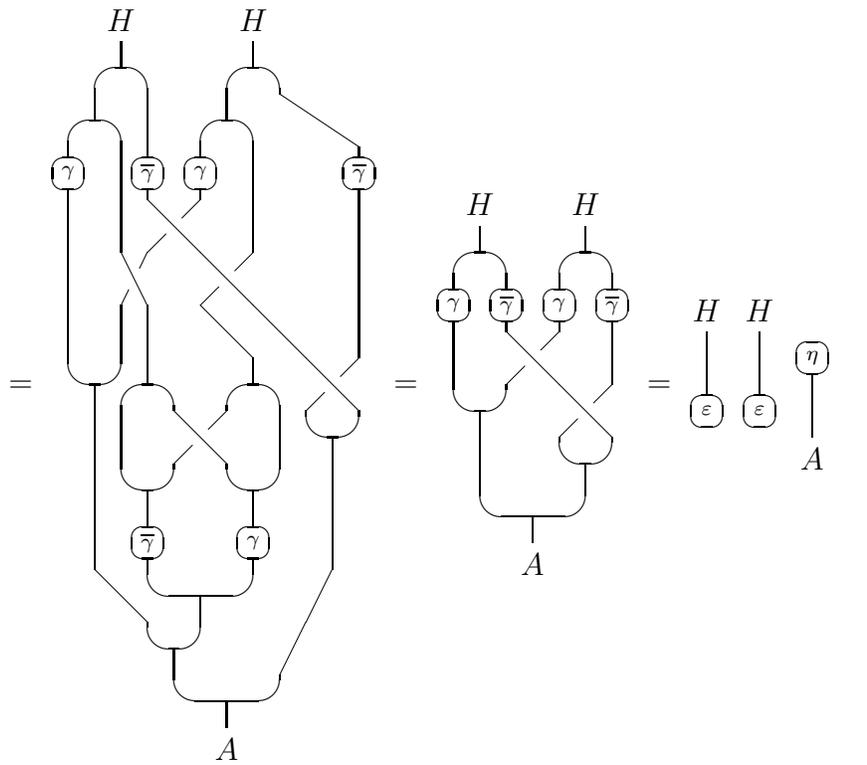
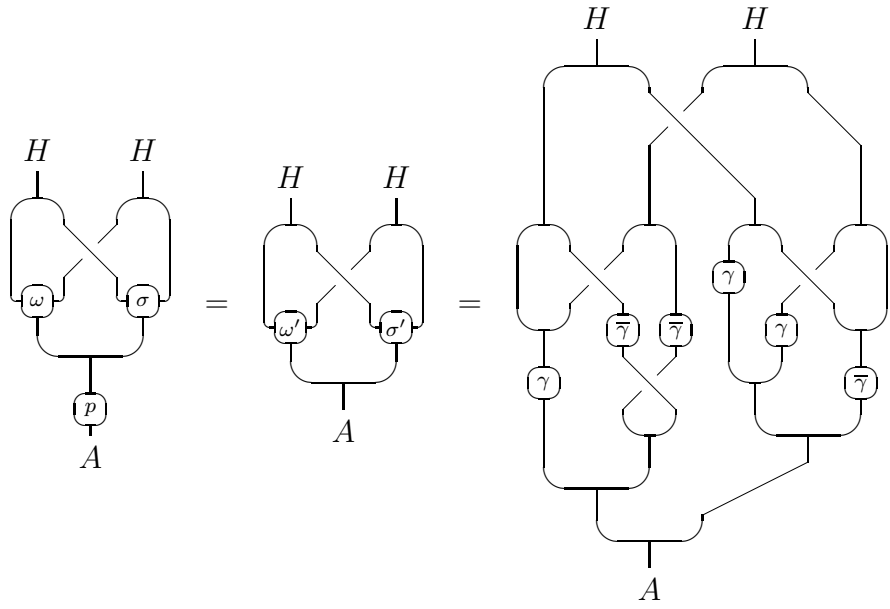
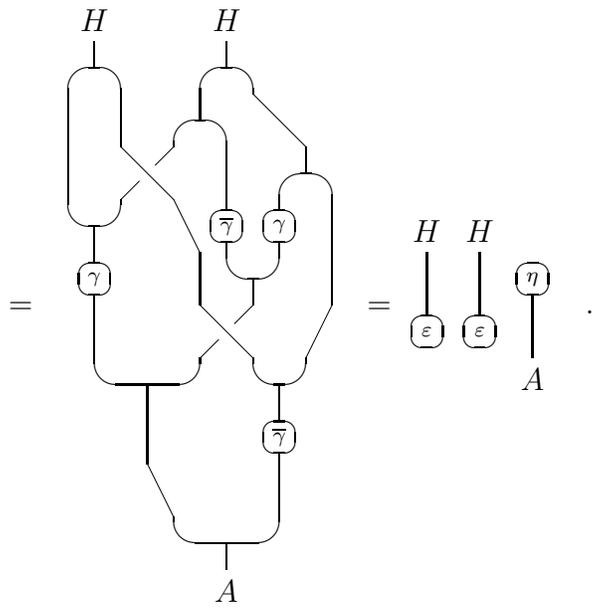
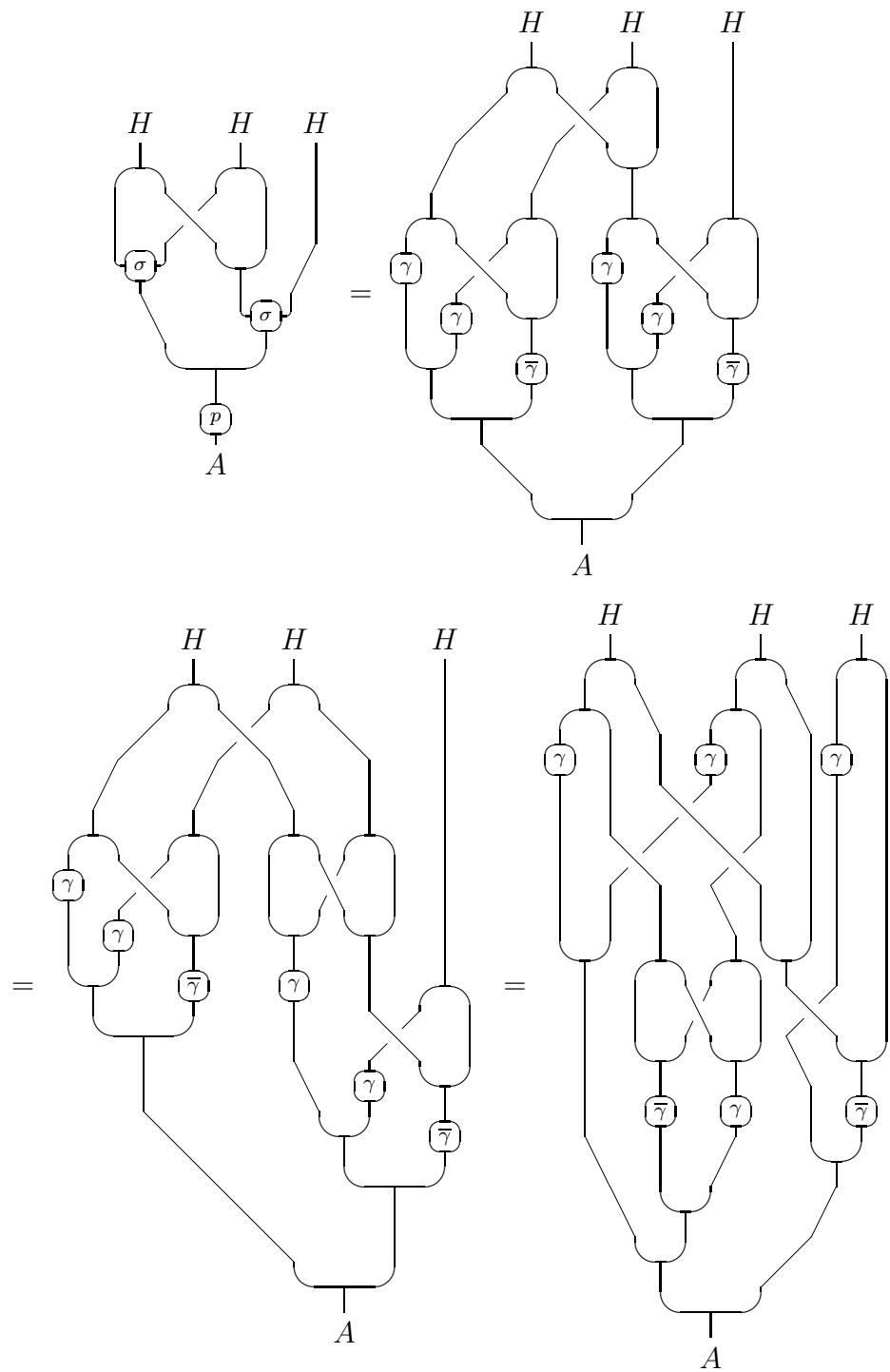
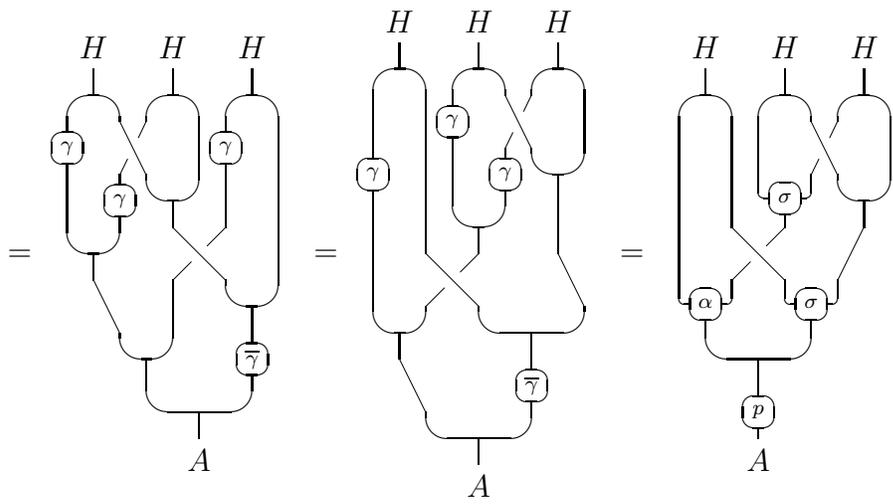
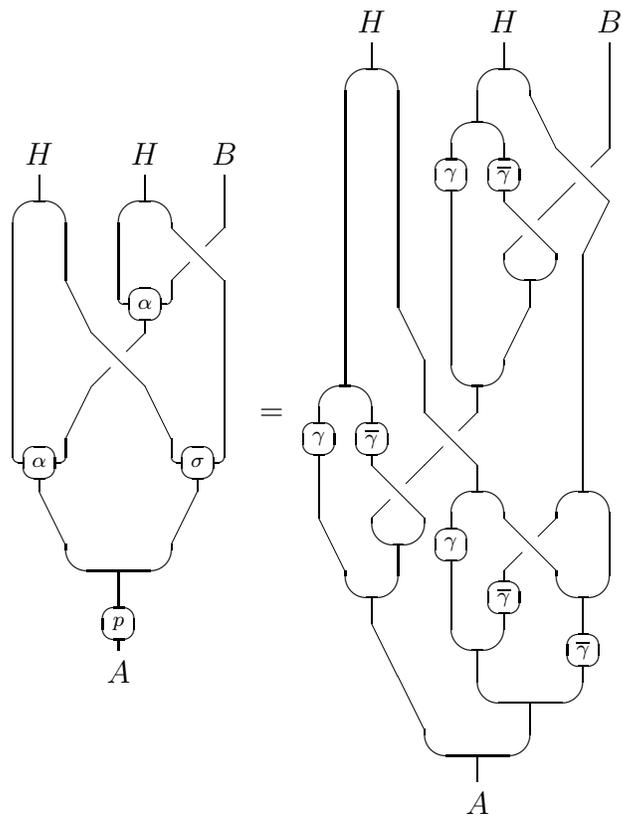
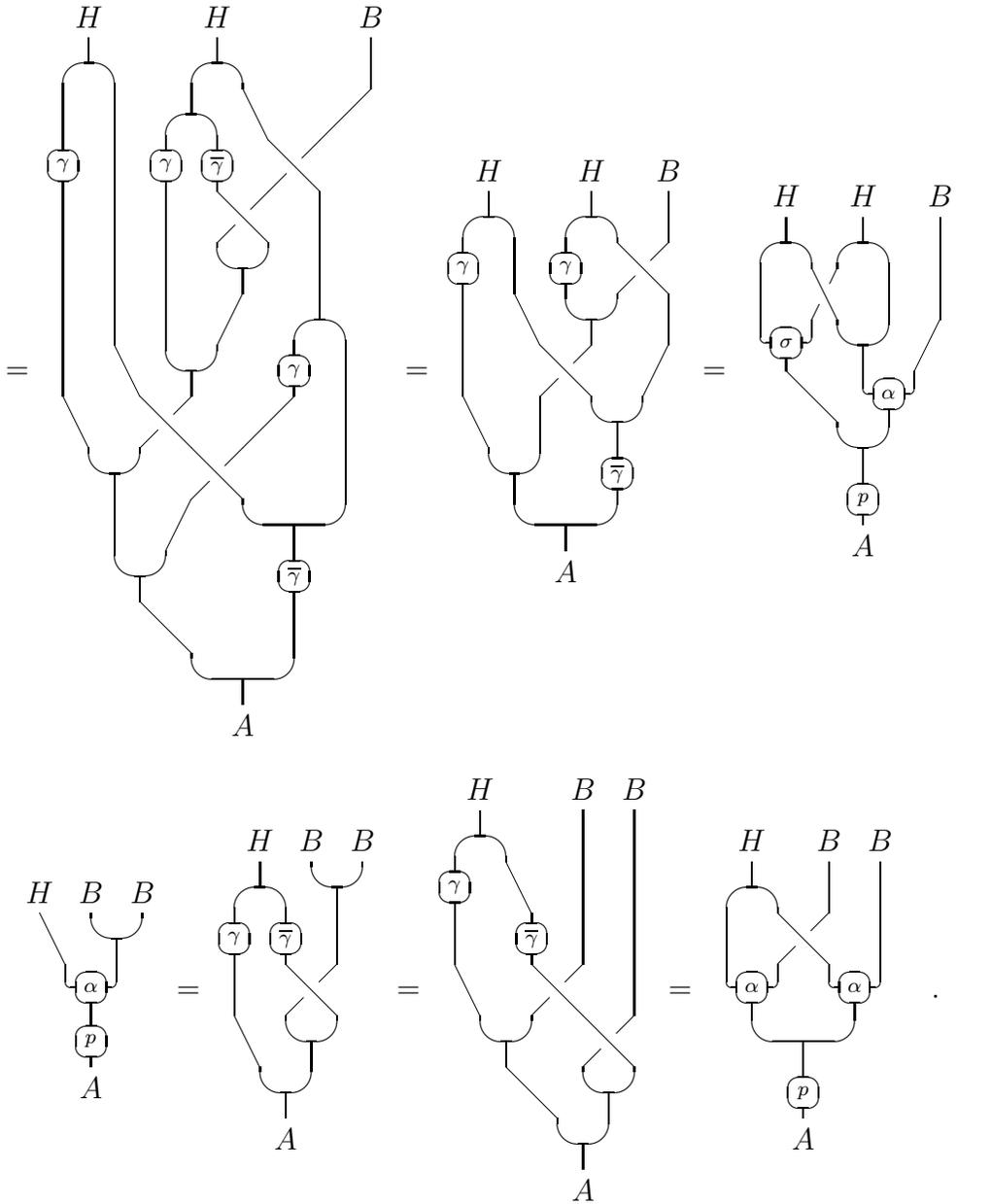
\[
\begin{tangle}
\object{H}\step[2]\object{H}\\
\tu {\sigma '}\\
\td \Theta\\
\object{A}\step[2]\object{H}\\
 \end{tangle}
  \step =\step
 \begin{tangle}
\step\object{H}\step[2]\step[2]\object{H}\\
\cd\step[2]\cd\\
\O \gamma\step[2]\x\step[2]\id\\
\id\step[2]\O \gamma\step[2]\cu\\
\cu\step[2]\step\O {\overline{\gamma} }\\
\step\Cu\\
\step[2]\td \Theta\\
\step[2]\object{A}\step[2]\object{H}\\
  \end{tangle}
  \step =\step
 \begin{tangle}
\step[2]\step[2]\object{H}\step[2]\step[2]\object{H}\\
\step[2]\step\cd\step[2]\cd\\
\step[2]\step\O \gamma \step[2]\x\step[2]\id\\
\step[2]\dd\step[2]\O \gamma\step[2]\cu\\
\step\dd\step[2]\step\id\step[2]\step\O {\overline{\gamma}}\\
\td \Theta\step[2]\td \Theta\step[2]\id\\
\id\step[2]\x\step[2]\id\step[2]\id\\
\cu\step[2]\cu\step\td \Theta\\
\step\id\step[4]\x\step[2]\id\\
\step[1]\Cu\step[2]\cu\\
\step[2]\step\object{A}\step[3]\step[2]\object{H}\\
 \end{tangle}
  \step =\step
 \begin{tangle}
\step[2]\step[2]\object{H}\step[2]\step[2]\object{H}\\
\step[2]\step\cd\step[2]\cd\\
\step[2]\step\O \gamma \step[2]\x\step[2]\id\\
\step[2]\dd\step[2]\O \gamma\step[2]\cu\\
\step\dd\step[2]\step\id\step[2]\cd\\
\td \Theta\step[2]\td \Theta\step\S\step[2]\O {\overline{\gamma}}\\
\id\step[2]\x\step[2]\id\step\id\step[2]\id\\\\
\cu\step[2]\cu\step\x\\
\step\id\step[4]\x\step[2]\id\\
\step[1]\Cu\step[1]\step\cu\\
\step[2]\step\object{A}\step[3]\step[2]\object{H}\\
 \end{tangle}\]
\vskip 0.1cm
\[
 \step =\step
 \begin{tangle}
\step[2]\step[2]\object{H}\step[2]\step[2]\object{H}\\
\step[2]\step\cd\step[2]\cd\\
\step[2]\step\id \step[2]\x\step[2]\id\\
\step[1]\step\ne2 \step[1]\ne1\step [2]\id\step[2]\nw1\\
\cd\step\cd\step\cd\step\cd\\
\id\step[2]\hx\step[2]\id\step\id\step[2]\hx\step[2]\id\\
\O \gamma\step[2]\O \gamma\step\cu\step\cu\step\cu\\
\cu\step[2]\id\step[2]\step\S\step\step[2]\O {\overline{\gamma}}\\
\step\id\step[2]\step\id\step[3]\id \step[2]\ne1\\
\step\id\step[2]\step\nw1\step[1]\step\x\\
\step\id\step[4]\x\step[2]\id\\
\step\Cu\step[2]\cu\\
\step[2]\step\object{A}\step[2]\step[3]\object{H}\\
  \end{tangle}
 \step =\step
 \begin{tangle}
\step[2]\step[2]\object{H}\step[2]\step[2]\object{H}\\
\step[2]\step\cd\step[2]\cd\\
\step[2]\step\id \step[2]\x\step[2]\id\\
\step[1]\step\ne2 \step[1]\ne1\step [2]\id\step[2]\nw1\\
\cd\step\cd\step\cd\step\cd\\ 
\id\step[2]\hx\step[2]\id\step\S\step[2]\hx\step[2]\id\\ 
\O \gamma\step[2]\O \gamma\step\cu\step\id\Step\S\step\cu\\ 
\cu\step[2]\id\step[2]\x\step[2]\O {\overline{\gamma}}\\ 
\step\id\step[2]\step\id\step\step\cu\step\ne1\\ 
\step\id\step[2]\step\nw1\step[1]\step\x\\
\step\id\step[4]\x\step[2]\id\\
\step\Cu\step[2]\cu\\
\step[2]\step\object{A}\step[2]\step[3]\object{H}\\
  \end{tangle}
   \step =\step
 \begin{tangle}
\step[2]\step[2]\object{H}\step[2]\step[2]\object{H}\\
\step[2]\step\cd\step[2]\cd\\
\step[2]\cd \step[1]\x\step[2]\id\\
\step \ne2 \step[2]\hx\step [2] \id \step [2]\nw1\\
\O \gamma\step[2]\step \ne2\step \id \step \cd\step[2]\id\\
\id\step[2]\O \gamma\step[2]\ne1\step[1]\S\step[2]\cu\\
\cu\step[2]\cu\step[2]\step\O {\overline{\gamma}}\\
\step\id\step[2]\step[2]\id \step [3]\ne2\\
\step\id\step[2]\step[2]\x\\
\step\Cu\step[2]\id\\
\step[2]\step\object{A}\step[2]\step[2]\object{H}\\
 \end{tangle}
 \]
\vskip 0.1cm
 \[
  \step =\step
 \begin{tangle}
\step\object{H}\step[2]\step[2]\object{H}\\
\cd\step[2]\cd\step[2]\Q \eta\\
\O \gamma\step[2]\x\step[2]\id\step[2]\id\\
\id\step[2]\O \gamma\step[2]\cu\step[2]\id\\
\cu\step[2]\step\O {\overline{\gamma} }\step[2]\step\id\\
\step\Cu\step[2]\step\id\\
\step[3]\object{A}\step[5]\object{H}\\
  \end{tangle}
    \step =\step
     \begin{tangle}
\object{H}\step[2]\object{H}\\
\tu {\sigma  '}\step \Q \eta\\
\step\object{A}\step[2]\object{H}\\
 \end{tangle}
 \Step . \step[6]
 \]
\vskip 0.1cm
 \[
\begin{tangle}
\object{H}\step[2]\object{B}\\
\tu {\alpha '}\\
\td \Theta\\
\object{A}\step[2]\object{H}\\
 \end{tangle}
  \step =\step
  \begin{tangle}
\step\object{H}\step[3]\object{B}\\
\cd\step[2]\id\\
\O \gamma\step[2]\O {\overline{\gamma}}\step[2]\O p\\
\id\step[2]\x\\
\d\step\cu\\
\step\cu\\
\step\td \Theta\\
\step\object{A}\step[2]\object{H}\\
 \end{tangle}
  \step =\step
  \begin{tangle}
  \step[2]\step\object{H}\step[3]\object{B}\\
  \step[2]\cd\step[2]\id\\
  \step[2]\O \gamma\step[2]\O {\overline{\gamma}}
  \step[2]\O p\\
  \step\dd\step[2]\x \step[2]\\
  \step\id \step[3]\id \step[2] \nw2 \step[2]\\
  \step\id\step[2] \td \psi\step[2] \td \psi\\
  \step\id\step[2]\id\step[2]\x\step[2]\id\\
 \td \psi\step\cu\step[2]\cu\\
  \id\step[2]\x\step[2]\sw2\\
  \cu\step[2]\cu\\
  \step\object{A}\step[2]\step[2]\object{H}\\
 \end{tangle} \step =\step
  \begin{tangle}
\step[2]\object{H}\step[3]\object{B}\\
\step\cd\step\step[1]\id\\
\step\O \gamma\step[2]\O {\overline{\gamma}}\step\step[1]\O p\\
\step\id\step[2]\x\step\\
\step\id\step[2]\id\step [2]\nw1\\
\step\id\step[2]\id\step[2] \td \psi\\
 \td \psi\step[]\cu\step\dd\\
\id\step[2]\x\step\step\id\\
\cu\step[2]\cu\\
  \step\object{A}\step[2]\step[2]\object{H}\\
 \end{tangle}
   \]
  \vskip 0.1cm

 \[
  \step =\step
  \begin{tangle}
\step[2]\object{H}\step[6]\object{B}\\
\step\cd\step[2]\step[2]\step\id\\
\step\id\step[2]\nw2\step[2]\step[2]\id\\
\cd\step[2]\cd\step[2]\O p\\
\O \gamma\step[2]\id\step[2]\S\step[2]\O {\overline{\gamma}}
\step[2]\id\\
\id\step[2]\cu\step[2]\x\\
\id\step[2]\step\nw2\step[2]\cu\\
\d\step[4]\x\\
\step\Cu\step[2]\id\\
\step[3]\object{A}\step[2]\step[2]\object{H}\\
 \end{tangle}
   \step =\step
  \begin{tangle}
\step\object{H}\step[3]\object{B}\\
\cd\step[2]\id\step[2]\Q \eta\\
\O \gamma\step[2]\O {\overline{\gamma}}\step[2]\O p\step[2]\id\\
\id\step[2]\x\step[2]\id\\
\d\step\cu\step[2]\id\\
\step\cu\step[2]\step\id\\
\step[2]\object{A}\step[4]\object{H}\\
 \end{tangle}
   \step =\step
  \begin{tangle}
\object{H}\step[2]\object{B}\\
\id\step[2]\id \Step\Q \eta\\
\tu \alpha\step[2]\id\\
\step\object{A}\step[3]\object{H}\step [2] .\\
 \end{tangle}
   \]
\vskip 0.1cm

\[
\begin{tangle}
\object{H}\step[2]\object{H}\\
\tu {\omega'}\\
\td \Theta\\
\object{A}\step[2]\object{H}\\
 \end{tangle}
  \step =\step
\begin{tangle}
\step\object{H}\step[3]\step\object{H}\\
\cd\step[2]\cd\\
\id\step[2]\x\step[2]\id\\
\cu\step[2]\O {\overline{\gamma}}\step[2]\O {\overline{\gamma}}\\
\step\O \gamma\step[2]\step\x\\
\step\id\step[2]\step\cu\\
\step\Cu\step\\
\step[2]\td \Theta\\
\step[2]\object{A}\step[2]\object{H}\\   
 \end{tangle}
  \step =\step
\begin{tangle}
\step[2]\object{H}\step[4]\step\object{H}\\
\step\cd\step[2]\step\cd\\

\step\id\step[2]\nw1\step[2]\id \step[2]\id\\

\step\id\step[3]\x\step[2]\id\\

\step\id\step[3]\id \step[2]\x \\

\cd\step\cd\step []\nw2 \step[1]\nw3\\

\id\step[2]\hx\step[2]\id\step[2]\cd\step\cd\\
\cu\step\cu\step[2]\S\step[2]\id\step\S\step[2]\O
{\overline{\gamma}}\\
\step\O \gamma\step[2]\step\id\step[2] \step\x\step\x\\
\step\id\step[2]\step\id\step[2]\step\O {\overline{\gamma}}
\step[2]\hx\step[2]\id\\
\step\id\step[2]\step\nw2\step[2]\cu\step\cu\\

\step\d\step[4]\x\step[2]\ne1\\

\step[2]\Cu\step[2]\cu\\
\step[2]\step\object{A}\step[2]\step[4]\object{H}\\
\end{tangle}
\]
\vskip 0.1cm
\[
   \step =\step
\begin{tangle}
\step\object{H}\step[3]\step\object{H}\\
\cd\step[2]\cd\step[2]\Q \eta\\
\id\step[2]\x\step[2]\id\step[2]\id\\
\cu\step[2]\O {\overline{\gamma}}\step[2]\O {\overline{\gamma}}\step[2]\id\\
\step\O \gamma\step[2]\step\x\step[2]\id\\
\step\id\step[2]\step\cu\step[2]\id\\
\step\Cu\step\step[2]\id\\
\step[2]\object{A}\step[2]\object{H}\\
 \end{tangle}
   \step =\step
  \begin{tangle}
\object{H}\step[2]\object{H}\\
\id\step[2]\id \Step\Q \eta\\
\tu \omega\step[2]\id\\
\step\object{A}\step[3]\object{H}\\
 \end{tangle}
 \step\step .\step[2]
  \]
\vskip 0.1cm
\[
 \begin{tangle}
\step\object{H}\step[4]\object{H}\\
\cd\step[2]\cd\\
\id\step[2]\x\step[2]\id\\
\tu \sigma\step[2]\tu \omega\\
\step\Cu\\
\step\step\step\O p\\
\step[2]\step\object{A}\\
  \end{tangle}
  \step =\step
   \begin{tangle}
\step\object{H}\step[4]\object{H}\\
\cd\step[2]\cd\\
\id\step[2]\x\step[2]\id\\
\tu {\sigma '} \step[2]\tu {\omega '}\\
\step\Cu\\
\step[2]\step\object{A}\\
  \end{tangle}
   \step =\step
 \begin{tangle}
\step[2]\step\object{H}\step[6]\object{H}\\
\step\Cd\step[2]\Cd\\
\step\id\step[2]\step[2]\x \step[2]\step[2]\nw2\\
\step\id\step[2]\step[2]\id \step[2]\nw2\step[5]\id\\
\cd\step[2]\cd\step[2]\cd\step[2]\cd\\
\O \gamma\step[2]\x\step[2]\id\step[2]\id\step[2]\x\step[2]\id\\
\id\step[2]\O \gamma\step[2]\cu\step[2] \cu\step[2]\O
{\overline{\gamma}}\step[2]\O {\overline{\gamma}}\\  
\cu\step[2]\step\O {\overline{\gamma} }\step\step[2] \step\O
\gamma\step[2]\step\x\\  
\step\Cu\step\step[2]\step\id\step[2]\step\cu\\
\step[2]\step\nw4\step[2]\step[2]\step\Cu\\
\step[7]\Cu\\
\step[9]\object{A}\\
  \end{tangle}
\]
\vskip 0.1cm
\[
 \step =\step
 \begin{tangle}
\step[2]\object{H}\step[5]\object{H}\\
\step\cd\step[2]\step\cd\\
\cd\step\id\step[2]\cd\step\nw3\\
\O \gamma\step[2]\id\step\O {\overline{\gamma}} \step[2] \O \gamma
\step[2]\id\step[2]\step[2]\O {\overline{\gamma}}\\
\id\step[2]\id\step\x\step[2]\id\step[2]\step[2]\id\\
\id \step[2]\hx\step[2]\x\step[4]\id\\
\id \step[2]\id\step[]\id\step[2]\nw2\step[] \nw2\step [3]\id\\
\cu\cd\step[2]\cd \step[1]\x\\
\step\id\step\id\step[2]\x\step[2]\id\step\cu\\
\step\id\step\cu\step[2]\cu\step[2]\id\\
\step\id\step[2]\O {\overline{\gamma}}\step[2]\step[2]
 \O \gamma\step[2]\step\id\\
 \step\nw2\step\Cu\step[2]\dd\\
 \step[2]\step\cu\step[2]\step\dd\\
 \step[2]\step[2]\Cu\\
\step[6]\object{A}\\
  \end{tangle}
 \step =\step
 \begin{tangle}
\step\object{H}\step[4]\object{H}\\
\cd\step[2]\cd\\
\O \gamma\step[2]\O {\overline{\gamma}} \step[2] \O \gamma
\step[2]\O {\overline{\gamma}}\\
\id\step[2]\x\step[2]\id\\
\cu\step[2]\x\\
\step\id\step[2]\step\cu\\
\step\Cu\\
\step[3]\object{A}\\
  \end{tangle}
 \step =\step
 \begin{tangle}
\object{H}\step[2]\object{H}\\
\id\step[2]\id\step[2]\Q \eta\\
\QQ \varepsilon\step[2]\QQ \varepsilon\step[2]\id\\
\step[2]\step[2]\object{A}\\
\end{tangle}
\step\step .
\]

\vskip 0.1cm
\[
 \begin{tangle}
\step\object{H}\step[4]\object{H}\\
\cd\step[2]\cd\\
\id\step[2]\x\step[2]\id\\
\tu \omega\step[2]\tu \sigma\\
\step\Cu\\
\step\step\step\O p\\
\step[2]\step\object{A}\\
  \end{tangle}
  \step =\step
   \begin{tangle}
\step\object{H}\step[4]\object{H}\\
\cd\step[2]\cd\\
\id\step[2]\x\step[2]\id\\
\tu {\omega '} \step[2]\tu {\sigma '}\\
\step\Cu\\
\step[2]\step\object{A}\\
  \end{tangle}
   \step =\step
 \begin{tangle}
\step[2]\step\object{H}\step[6]\object{H}\\
\step\Cd\step[2]\Cd\\
\step\id\step[2]\step[2]\x \step[2]\step[2]\nw2\\
\step\id\step[2]\step[2]\id \step[2]\nw2\step[5]\id\\
\cd\step[2]\cd\step[2]\cd\step[2]\cd\\
\id\step[2]\x\step[2]\id\step[2]\O \gamma\step[2]\x\step[2]\id\\
\cu\step[2]\O {\overline{\gamma}}\step[2]\O {\overline{\gamma}}
  \step[2]   \id\step[2]\O \gamma\step[2]\cu \\  
\step\O \gamma\step[2]\step\x\step[2]
\cu\step[2]\step\O {\overline{\gamma} }\\  
\step\id\step[2]\step\cu\step[2] \step\Cu\\
\step\Cu\step[2]\step[2]\step\ne4\\
\step[3]\Cu\\
\step[5]\object{A}\\
  \end{tangle}
\]
\vskip 0.1cm
\[
   \step =\step
 \begin{tangle}
\step\object{H}\step[2]\step[3]\object{H}\\
\cd\step[2]\step\cd\\
\id\step[2]\id\step[2]\cd\step\nw2\\
\id\step[2]\x\step[2]\id\step[2]\cd\\
\cu\step[2]\d\step[]\O {\overline{\gamma}} \step[2]\O \gamma
\step[2]\id\\
\step\O \gamma\step[2]\step[2]\id \step[]\cu\step[2]\id\\
\step\id \step[4]\x\step[2]\ne1\\
\step[]\Cu\step[2]\cu\\
\step[2]\step\id\step[2]\step[2]\step[]\O {\overline{\gamma}}\\
\step[2]\step\d\step[2]\step[2]\id\\
\step[2]\step[2]\Cu\\
\step[3]\step[2]\step\object{A}\\
\end{tangle}
   \step =\step
 \begin{tangle}
\object{H}\step[2]\object{H}\\
\id\step[2]\id\step[2]\Q \eta\\
\QQ \varepsilon\step[2]\QQ \varepsilon\step[2]\id\\
\step[2]\step[2]\object{A}\\
\end{tangle}
\step\step .
\]
\vskip 0.1cm
\[
 \begin{tangle}
\step\object{H}\step[2]\step[2]\object{H}\step\step[2]\object{H}\\
\cd\step[2]\cd\step[2]\id\\
\id\step[2]\x\step[2]\id\step[2]\id\\
\tu \sigma\step[2]\cu\step\dd\\
\step\d\step[2]\step\tu \sigma\\
\step[2]\Cu\\
\step[2]\step[2]\O p\\
\step[2]\step[2]\object{A}\\
\end{tangle}
\step =\step
 \begin{tangle}
\step[2]\step[2]\step\object{H}\step[2]\step[2]\object{H}\step\step[3]\object{H}\\
\step[2]\step[2]\cd\step[2]\cd\step[2]\step\id\\
\step[2]\sw2\step[2]\step\x\step[2]\id\step[2]\step\id\\
\step\dd\step[2]\step\dd\step[2]\cu\step[2]\step\id\\
\cd\step[2]\cd \step[2]\cd\step[2]\cd\\   
\O \gamma\step[2]\x\step[2]\id\step[2]\O \gamma\step[2]\x\step[2]\id\\
\id\step[2]\O \gamma\step[2]\cu\step[2]\id\step[2]\O \gamma\step[2]\cu\\
\cu\step[2]\step\O {\overline{\gamma}}\step\step[2]
\cu\step[2]\step\O {\overline{\gamma}}\step\\ 
\step\Cu\step \step[2] \step\Cu\step\\  
\step[2]\step\nw2\step[2]\step[2]\step\sw2\\
\step[2]\step[2]\step\Cu\\
\step[3]\step[4]\object{A}\\
\end{tangle}
\]
\vskip 0.1cm
\[
\step =\step
 \begin{tangle}
\step[2]\step[2]\step\object{H}\step[2]\step[2]\object{H}\step\step[5]\object{H}\\
\step[2]\step[2]\cd\step[2]\cd\step[2]\step[2]\step\id\\
\step[2]\sw2\step[2]\step\x\step[2]\nw2\step[2]\step[2]\id\\
\step\dd\step[2]\step\dd\step[2]\d\step[2]\step\id\step[2]\step\id\\
\cd\step[2]\cd \step[2]\cd\step\cd\step[2]\id\\   
\O \gamma\step[2]\x\step[2]\id\step[2]\id\step[2]\hx\step[2]\id\step[2]\id\\
\id\step[2]\O \gamma\step[2]\cu\step[2]\cu\step\cu\step[2]\id\\
\cu\step[2]\step\O {\overline{\gamma}}\step\step[2]
\step\O \gamma\step\step[2]\id\step[2]\cd\\ 
\step\Cu\step \step[2] \step\id\step\step[2]\x\step[2]\id\\  
\step[2]\step\id\step[2]\step[2]\step[2]\d \step[2]\O \gamma
\step[2]\cu\\
\step[2]\step\nw2\step\step[5]\cu\step[2]\step\O {\overline{\gamma}}\\
\step[5]\nw2\step\step[2]\step[2]\Cu\\
\step[7]\nw2\step\step[2]\step[2]\id\\
\step[9]\Cu\\
\step[6]\step[5]\object{A}\\
\end{tangle}
\step =\step
 \begin{tangle}
\step[2]\object{H}\step[2]\step[4]\object{H}\step\step[3]\object{H}\\
\step\cd\step[2]\step[2]\cd\step[2]\cd\\
\cd\step\d\step[2]\cd\step\d\step\id\step[2]\id\\
\O \gamma\step[2]\id\step[2]\id\step[2]\O
\gamma\step[2]\id\step[2]\id\step\O \gamma\step\step[1]\id\\
\id\step[2]\id \step[2]\x\step[2]\id\step[2]\id\step\id\step[2]\id\\
\id\step[2]\x\step\step[1]\x\step[2]\id\step\id\step[2]\id\\
\id\step[2]\id \step\step[1]\id\step[2]\nw1 \step\id\step[2]\id\step\id\step[2]\id\\
\cu\step\cd\step\cd\cu\step\id\step[2]\id\\
\step\id\step[2]\id\step[2]\hx\step[2]\id\step\x\step[2]\id\\
\step\id\step[2]\cu\step\cu\step[1]\nw1\step\cu\\
\step\id\step[2]\step\O {\overline{\gamma}}\step[2]\step \O \gamma
\step[2]\step[1]\id\step[2] \O {\overline{\gamma}}\\
\step\d\step[2]\id\step[2]\dd\step[2]\step[1]\cu\\
\step[2]\d\step\cu\step[2]\step[2]\ne1\\
\step[2]\step\cu\step[2]\step\sw2\\
\step[2]\step[2]\Cu\\
\step[6]\object{A}\\
\end{tangle}
\]

\vskip 0.1cm
\[
\step =\step
 \begin{tangle}
\step\object{H}\step[2]\step\object{H}\step[3]\object{H}\\
\cd\step\cd\step\cd\\
\O \gamma\step[2]\hx\step[2]\id\step\O \gamma\step[2]\id\\
\id\step[2]\O \gamma\step\cu\step\id\step[2]\id\\
\cu\step[2]\x\step[2]\id\\
\step\d\step[2]\id\step[2]\cu\\
\step[2]\cu\step[2]\step\O {\overline{\gamma}}\\
\step\step[2]\Cu\\
\step[5]\object{A}\\
\end{tangle}
\step =\step
 \begin{tangle}
\step\object{H}\step[2]\step\object{H}\step[3]\object{H}\\
\cd\step\cd\step\cd\\
\id\step[2]\id\step\O \gamma\step[2]\hx\step[2]\id\\
\O \gamma\step[2]\id\step\id\step[2]\O \gamma\step\cu\\
\id\step[2]\id\step\cu\step[2]\id\\
\id\step[2]\x\step[2]\step\nw1\\
\cu\step[]\step\Cu\\
\step\d\step[2]\step[2]\O {\overline{\gamma}}\\
\step[2]\Cu\\
\step[4]\object{A}\\
\end{tangle}
\step =\step
 \begin{tangle}
\step\object{H}\step[2]\step[]\object{H}\step[3]\object{H}\\
\cd\step[]\cd\step\cd\\
\id\step[2]\id\step[1]\id\step[2]\hx\step[2]\id\\
\id\step[2]\id\step[1]\tu \sigma\step\cu\\
\id\step[2]\x\step[2]\ne1\\
\tu \alpha\step[2]\tu \sigma\\
\step\Cu\\
\step[2]\step\O p\\
\step[3]\object{A}\\
\end{tangle}
\step[2]\step .\step[4]\step[2]
\]
\vskip 0.1cm
\[
 \begin{tangle}
\step\object{H}\step[3]\step\object{H}\step[3]\object{B}\\
\cd\step\step\cd\step\step\id\\
\id\step[2]\id\step\step\id\step[2]\x\\
\id\step[2]\nw1\step\tu \alpha\step\step\id\\
\id\step[3]\x\step\step[2]\id\\
\id\step[2]\ne1 \step[2]\nw1\step[2]\id\\
\tu \alpha\step[2]\step[2]\tu \sigma\\
\step\d\step[2]\step[2]\dd\\
\step[2]\Cu\\
\step[2]\step[2]\O p\\
\step[4]\object{A}\\
\end{tangle}
\step =\step
 \begin{tangle}
\step\step\object{H}\step[4]\step\object{H}\step[4]\object{B}\\
\step\cd\step[2]\step\cd\step[2]\step\id\\
\step\id\step[2]\id\step[2]\cd\step\d\step[2]\id\\
\step\id\step[2]\id\step[2]\O \gamma\step[2]\O {\overline{\gamma}}
\step[2]\x\\  
\step\id\step[2]\id\step[2]\id\step[2]\x\step[1]\ne1\\
\step\id\step[2]\id\step[2]\id\step[2]\cu\step\id\\
\step\id\step[2]\d\step\id\step[2]\dd\step[2]\id\\  
\cd\step[2]\id\step \cu\step[2]\step\id\\
\O \gamma\step[2]\O
{\overline{\gamma}}\step[2]\x\step[2]\step[2]\id\\  
\id\step[2]\x\step\cd\step[2]\cd\\
\d\step\cu\step\O \gamma\step[2]\x\step[2]\id\\
\step\cu\step[2]\id\step[2]\O {\overline{\gamma}}\step[2]\cu\\  
\step[2]\d\step[2]\cu\step[2]\step\O {\overline{\gamma}}\\
\step[2]\step \d\step[2]\Cu\\
\step[2]\step[2]\Cu\\
\step[6]\object{A}\\
\end{tangle}
\]
\vskip 0.1cm
\[
\step =\step
 \begin{tangle}
\step\object{H}\step[5]\object{H}\step[6]\object{B}\\
\cd\step[2]\step\cd\step[2]\step[2]\step\id\\
\id\step[2]\id\step[2]\cd\step\d\step[2]\sw2\\
\O \gamma\step[2]\id\step[2]\O \gamma\step[2]\O
{\overline{\gamma}}\step[2]\x\\
\id\step[2]\id\step[2]\id\step[2]\x\step[2]\id\\  
\id\step[2]\id\step[2]\id\step[2]\cu\step[2]\id\\
\id\step[2]\id\step[2]\id\step[2]\dd\step[2]\cd\\
\id\step[2]\nw1\step\cu\step[2]\step[1]\O \gamma\step[2]\id\\  
\nw1\step[2]\x\step\step\sw2\step[2]\step\id\\
\step \cu\step[2]\x\step[2]\step[2]\id\\  
\step[2]\id \step[2]\ne1\step[1]\step\Cu\\
\step[2]\cu\step[2]\step[2]\step\O {\overline{\gamma}}\\
\step[2]\step\nw2\step[2]\step[2]\step\id\\
\step[5]\Cu\\
\step[7]\object{A}\\
\end{tangle}
\step =\step
 \begin{tangle}
\step\object{H}\step[3]\step\object{H}\step[3]\object{B}\\
\cd\step\step\cd\step\step\id\\
\O \gamma\step[2]\id\step\step\O \gamma\step[2]\x\\
\id\step[2]\nw1\step\cu\step\step\id\\
\id\step[3]\x\step\step[]\ne1\\
\nw1 \step[2]\id\step[2]\cu\\
\step \cu \step\step[2]\O {\overline{\gamma}}\\
\step[2]\Cu\\
\step[4]\object{A}\\
\end{tangle}
\step =\step
 \begin{tangle}
\step\object{H}\step[3]\object{H}\step[3]\object{B}\\
\cd\step\cd\step[2]\id\\
\id\step[2]\hx\step[2]\id\step[2]\id\\
\tu \sigma\step\cu\step\dd\\
\step\nw2\step[2]\tu \alpha\\
\step[2]\step\cu\\
\step[2]\step[2]\O p\\
\step[4]\object{A}\\
\end{tangle}
\step\step .
\]

\vskip 0.1cm
\[
\begin{tangle}
\object{H}\step[2]\object{B}\step[2]\object{B}\\
\d\step\cu\\
\step\tu \alpha\\
\step[2]\O p\\
\step[2]\object{A}\\
\end{tangle}
\step =\step
\begin{tangle}
\step\object{H}\step[2]\object{B}\step[2]\object{B}\\
\cd\step\cu\\
\O \gamma\step[2]\O {\overline{\gamma}}\step[2]\id\\
\id\step[2]\x\\
\d\step\cu\\
\step\cu\\
\step[2]\object{A}\\
\end{tangle}
\step =\step
\begin{tangle}
\step\object{H}\step[2]\step[2]\object{B}\step[2]\object{B}\\
\cd\step[2]\step\id\step[2]\id\\
\O \gamma\step[2]\d\step[2]\id\step[2]\id\\
\id\step[2]\step\O {\overline{\gamma}}\step[2]\id\step[2]\id\\
\d\step[2]\x\step[2]\id\\
\step\cu\step[2]\x\\
\step[2]\nw2\step\step\cu\\
\step[2]\step[2]\cu\\
\step[5]\object{A}\\
\end{tangle}
\step =\step
\begin{tangle}
\step\object{H}\step[2]\step\object{B}\step[2]\object{B}\\
\cd\step[2]\id\step[2]\id\\
\id\step[2]\x\step[2]\id\\
\tu \alpha\step[2]\tu \alpha\\
\step\Cu\\
\step[2]\step\O p\\
\step[3]\object{A}\\
\end{tangle}
\step\step .
\]

\vskip 0.1cm \noindent Thus (2-COC) , (TM) and (WA) hold.
 Finally
,we show that $\Phi$ is an algebra isomorphism from $B \# _\sigma
H$ onto $A$ . \vskip 0.1cm
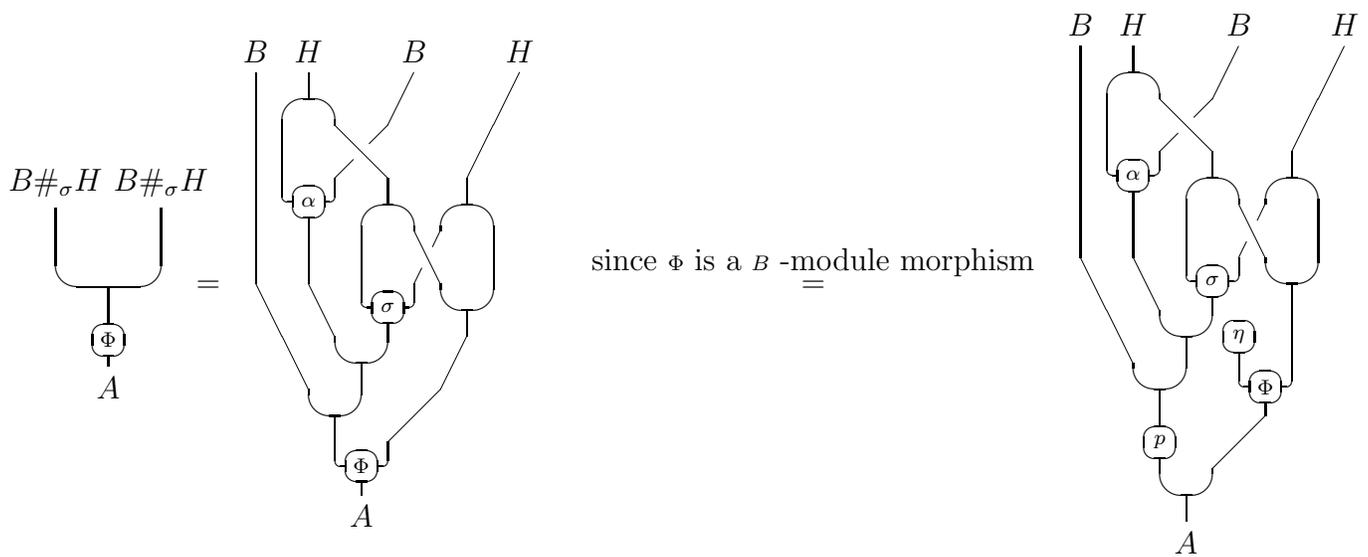
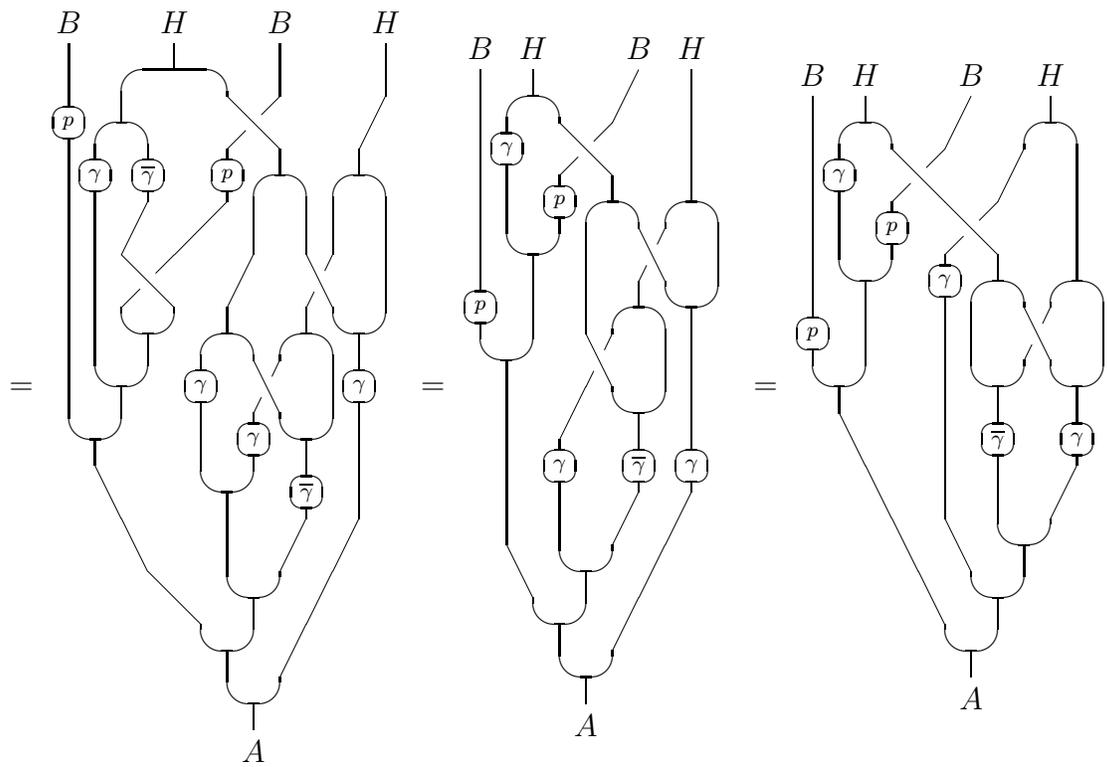
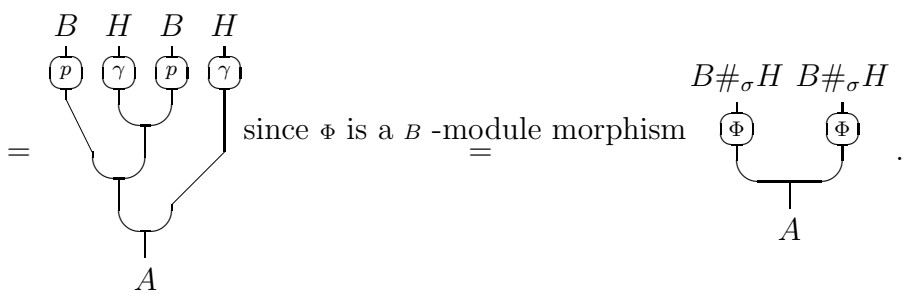
\[
\begin{tangle}
\object{B \# _\sigma H}\step[2]\step[2]\object{B \# _\sigma H}\\
\id\step\step\step\step\id\\
\Cu\\
\step\step\O \Phi\\
\step\step\object{A}\\
\end{tangle}
\step =\step
\begin{tangle}
\object{B}\step[2]\object{H}\step[2]\step[2]\object{B}\step[4]\object{H}\\
\id\step\cd\step[2]\dd\step[2]\step\dd\\
\id\step\id\step[2]\x\step[2]\step\dd\\
\id\step\tu \alpha\step\cd\step\cd\\
\id\step[2]\id\step[2]\id\step[2]\hx\step[2]\id\\
\d\step\d\step\tu \sigma\step\cu\\
\step\d\step\cu\step[2]\dd\\
\step[2]\cu\step\sw2\\
\step\step[2]\tu \Phi\\
\step[4]\object{A}\\
\end{tangle}
\step[1] \step \stackrel { \hbox {  since  } \Phi  \hbox { is a }
B \hbox { -module morphism }}{ =}\step
\begin{tangle}
\object{B}\step[2]\object{H}\step[2]\step[2]\object{B}\step[4]\object{H}\\
\id\step\cd\step[2]\dd\step[2]\step\dd\\
\id\step\id\step[2]\x\step[2]\step\dd\\
\id\step\tu \alpha\step\cd\step\cd\\
\id\step[2]\id\step[2]\id\step[2]\hx\step[2]\id\\
\d\step\d\step\tu \sigma\step\cu\\
\step\d\step\cu\step\Q \eta\step[2]\id\\
\step[2]\cu\step[2]\tu \Phi\\
\step\step[2]\O p\step[2]\sw2\\
\step[2]\step\cu\\
\step[4]\object{A}\\
\end{tangle}
\]
\vskip 0.1cm
\[
\step =\step
\begin{tangle}
\object{B}\step[4]\object{H}\step[2]\step[2]\object{B}\step[4]\object{H}\\
\id\step[2]\Cd\step[2]\id\step[2]\step[2]\id\\
\O p\step\cd\step[3]\x \step[2]\step\ne1\\
\id\step\O \gamma\step[2]\O {\overline{\gamma}}\step[3]\O
p\step[1]\cd\step\cd\\   
\id\step\id\step[]\ne1 \step[2]\ne2\step[] \id \step[2] \id \step \id  \step[2]\id\\
\id\step\id\step[]\x\step[2]\dd\step[2]\hx\step[2]\id\\
\id\step\id\step[]\cu\step\cd\step\cd\cu\\  
\id\step\cu\step[2]\O \gamma\step[2]\hx\step[2]\id\step\O \gamma\\
\cu\step[2]\step\id\step[2]\O \gamma\step\cu\step\id\\   
\step\d\step[2]\step\cu\step[2]\O {\overline{\gamma}}\step[2]\id\\  
\step[2]\d\step[2]\step\id\step[2]\dd\step\dd\\
\step[2]\step\nw2\step[2]\cu\step\dd\\  
\step[5]\cu\step\dd\\
\step[6]\cu\\
\step[7]\object{A}\\
\end{tangle}
\step =\step
\begin{tangle}
\object{B}\step[2]\object{H}\step[4]\object{B}\step[2]\object{H}\\
\id\step[1]\cd\step[2]\dd\step[2]\id\\
\id \step[1]\O \gamma\step[2]\x\step[2]\step\id\\
\id\step[1]\id\step[2]\O p\step\cd\step\cd\\
\id\step[1]\cu\step\id\step[2]\hx\step[2]\id\\
\O p \step[2]\id \step [2]\id\step\cd\cu\\
\cu \step[2]\hx\step[2]\id\step\id\\
\step\id\step[2]\dd\step\cu\step\id\\  
\step\id\step[2]\O \gamma\step[2]\step\O
{\overline{\gamma}}\step[2]\O \gamma\\
\step\id\step[2]\id\step[2]\dd\step\dd\\
\step\d\step\cu\step\dd\\
\step[2]\cu\step\dd\\
\step[2]\step\cu\\
\step[4]\object{A}\\
\end{tangle}
\step =\step
\begin{tangle}
\object{B}\step[2]\object{H}\step[4]\object{B}\step[3]\object{H}\\
\id\step\cd\step[2]\dd\step[2]\cd\\
\id\step\O \gamma\step[2]\x\step[2]\dd\step[2]\id\\
\id\step\id\step[2]\O p\step[2]\x\step[2]\step\id\\
\id\step\cu\step[2]\O \gamma\step\cd\step\cd\\
\O p \step [2] \id \step[2]\step\id\step\id\step[2]\hx\step[2]\id\\
\cu \step \step[2]\id\step\cu\step\cu\\
\step\d\step[2]\step\id\step[2]\O
{\overline{\gamma}}\step\step[2]\O
\gamma\\
\step[2]\d\step[2]\id\step[2]\id\step[2]\dd\\
\step[2]\step\d\step\d\step\cu\\
\step[2]\step[2]\d\step\cu\\
\step[5]\cu\\
\step[6]\object{A}\\
\end{tangle}
\]
\vskip 0.1cm
\[ \step = \step
\begin{tangle}
\object{B}\step[2]\object{H}\step[2]\object{B}\step[2]\object{H}\\
\O p\step[2]\O \gamma\step[2]\O p\step[2]\O \gamma\\
\d\step\cu\step[2]\id\\
\step\cu\step\sw2\\
\step[2]\cu\\
\step[3]\object{A}\\
\end{tangle}
\stackrel { \hbox {  since  } \Phi  \hbox { is a } B \hbox {
-module morphism } \ \ \ }{ =}
\begin{tangle}
\object{B \# _\sigma H}\step[2]\step[2]\object{B \# _\sigma H}\\
\O \Phi\step\step\step\step\O \Phi\\
\Cu\\
\step\step\object{A}\\
\end{tangle}
\step[2] . \]

Since $\Phi$ is an isomorphism we have  $\Phi (\eta _B \otimes
\eta _H) = \eta _A.$ \ \
 \begin{picture}(5,5)
\put(0,0){\line(0,1){5}}\put(5,5){\line(0,-1){5}}
\put(0,0){\line(1,0){5}}\put(5,5){\line(-1,0){5}}
\end{picture}

\section {Hopf  Galois extension  in the Yetter-Drinfeld  category}
Using the  conclusion in Theorem \ref {1.3}, in this section  we
give the relation between crossed product and $H$-Galois extension
in the Yetter-Drinfeld category $ {}^D_D {\cal YD}$ with Hopf
algebra $D$ over field $k$.

It follows from \cite [Corollary 2.2.8]{DNR01} that $^D_D {\cal
YD}$ is an additive category. Therefore , by  \cite [Page
242]{Fa73},  $equivaliser (f_1, f_2) = ker (f_1 - f_2)$ and
 $coequivaliser (f_1,  f_2) = A / Im (f _1 -f_2)$ for any two
 morphisms $f_1$, $f_2$ in  $^D_D {\cal YD}$.
Consequently we have

\begin {Corollary}\label {2.1}

Let ${\cal C}$ be the Yetter-Drinfeld category $^D_D {\cal YD}$.
Then the following assertions are equivalent:

(i) There exists an invertible 2-cocycle $\sigma : H \otimes H
\rightarrow B, $ and a week action $\alpha $ of $H$ on $B$ such
that $A = B \# _\sigma H$ is a crossed product algebra.

(ii) $(A, \psi )$ is a right $H$-comodule algebra with $B=A ^{co
H}$, the extension $A/B$ is  Galois, the canonical epic $q:
A\otimes A \rightarrow A\otimes_B A$ is split and $A$ is
isomorphic as left $B$-modules and right $H$-comodules to
$B\otimes H$ in ${\cal C}$, where the $H$-comodule operation and
$B$-module operation of $B \otimes H$ are   $\psi _{B \otimes H} =
id _B \otimes  \Delta _H$ and $\alpha _{B \otimes H} = m_B \otimes
id _H$ respectively.
\end {Corollary}
Corollary \ref {2.1} implies \cite [Theorem 6.4.12] {DNR01} since
$q$ is split in the category of vector spaces with trivial
braiding.

\begin {Example}\label {2.3}
Let $H$ be  a braided Hopf algebra in ${}^D_D{\cal YD}$.
  Let $B =H$ and   $H$ act  on $B$ by adjoint action
$\alpha = m (H \otimes m) (H \otimes C_{H, B})(H \otimes S \otimes
B) (\Delta _H \otimes B)$. Thus $B \# H=A$ is a smash product
algebra in  ${}^D_D{\cal YD}$. By Corollary \ref {2.1}, $A/B$ is
an $H$-Galois extension in  ${}^D_D{\cal YD}$.
\end {Example}

Note that many braided Hopf algebras have been known. For example,
 by the Radford's method in \cite [Theorem 1 and
Theorem 3]{Ra85}, one can obtain a braided Hopf algebra $H$ in the
Yetter-Drinfeld category ${}^D_D{\cal YD}$ for any graded Hopf
algebra $A$ with $D= A_0.$  One can also obtain a braided Hopf
algebra $\underline H$, the braided group analogue of $H$ in the
Yetter-Drinfeld category ${}^D_D{\cal YD}$ for any
(co)quasitriangular  Hopf algebra $H$ with $D= H$ (see \cite
{Ma95b}). Furthermore, there exist many graded Hopf algebras as they
can be constructed by Hopf quivers (see \cite {CR02}. Otherwise,
Nichols algebras also are braided Hopf algebrs.

\vskip 0.5cm {\bf Acknowledgement }: The first two authors were
financially supported by the Australian Research Council. S.C.Z
thanks the Department of Mathematics, University of Queensland and
Hong Kong University of Science and Technology for hospitality.
Thank Y.Bespalov and V.Lyubashenko for their t-angles.sty.


\providecommand{\bysame}{\leavevmode\hbox
to3em{\hrulefill}\thinspace}
\begin{thebibliography}{1}

\bibitem{Bes:cross}
Yu.~N. Bespalov, \emph{Crossed modules and quantum groups in
braided
  categories}, Applied Categorical Structures \textbf{5} (1997), no.~2,
  155--204, http://xxx.lanl.gov/abs/q-alg/9510013.

\bibitem{Bes:FRT}
\bysame, \emph{On the braided {FRT}-construction}, J. Nonlinear
Math. Phys.
  \textbf{4} (1997), no.~1-2, 195--205, http://xxx.lanl.gov/abs/q-alg/9510012.

\bibitem{BKLT:int}
Yu.~N. Bespalov, T.~Kerler, V.~V. Lyubashenko, and V.~G. Turaev,
  \emph{Integrals for braided {H}opf algebras}, J. Pure and Appl. Algebra
  \textbf{148} (2000), no.~2, 113--164, Available as
  http://xxx.lanl.gov/abs/q-alg/9709020.

\end{thebibliography}


\begin{thebibliography}{150}


 \bibitem {AS02} N. Andruskewisch and H.J.Schneider, Pointed Hopf algebras,
new directions in Hopf algebras, edited by S. Montgomery and H.J.
Schneider, Cambradge University Press, 2002.
\bibitem {AS00} N. Andruskewisch and H.J.Schneider, Finite quantum groups and Cartan matrices,
Advances in Mathematics, 154(2000), 1-45.
 \bibitem {BD99} Y. Bespalov, B. Drabant,  Cross product bialgebras - Part II,
   J. Aal. {\bf 240} (2001), 445-504.


\bibitem {BD98} Y.Bespalov, B.Drabant, Hopf  (bi-)modules and crossed modules
in braided monoidal categories, J. Pure and Applied Alg. {\bf 123}
(1998), 105-129.

\bibitem {BM89}  R. J. Blattner and   S. Montgomery, Crossed products and Galois extensions of Hopf algebras,
 Pacific J. Math. {\bf 137} (1989), 37-54.

\bibitem {CHR65} S.U.Chase, D.K.Harrison, and A.Rosenberg, Galois theory and cohomology of commutative
 rings, AMS Memoirs No.52, 1965.
\bibitem {CS69} S.U.Chase and M.E.Sweedler, Hopf algebras and Galois theory,
Lecture notices in Math 97, Springer, Berlin, 1969.

\bibitem {CR02} C.  Cibils, M.  Rosso,  Hopf quivers,  J. Alg. {\bf 254}
(2002), 241-251.

\bibitem {DNR01} S.Dascalescu, C.Nastasecu and S. Raianu,
Hopf algebras: an introduction,  Marcel Dekker Inc. , 2001.

\bibitem {DT89} Y.Doi and M.Takeuchi, Hopf-Galois extensions of algebras,
the Miyashita-Ulbrich action, and Azumaya algebras, J. Alg.
{\bf 121} (1989), 488-516.

\bibitem {DZ99} A.Van Daele and Y.Zhang, Galois theory for multiplier Hopf
algebras with integrals, Alg. and Rep. Theory, {\bf 2}
(1999), 83-106.

\bibitem {Fa73} C.Faith, Algebra: rings, modules and categries,
Springer-Verlag, 1973.

\bibitem {Ke99}  T.Kerler, Bridged links and tangle presentations of
cobordism categories, Adv. Math. {\bf 141} (1999), 207-- 281.

\bibitem {KT81}  H.F.Kreimer and M.Takeuchi, Hopf algebras and Galois
extensions of an algebras.
Indiana University Math. J. {\bf 30} (1981), 675--692.

\bibitem {Ma90a}  S. Majid,
Physics for algebraists: Non-commutative and non-cocommutative
Hopf algebras by a bicrossproduct construction, J. Alg.   {\bf
130}  (1990), 17--64.
 \bibitem {Ma95a} S. Majid. Algebras and Hopf algebras
  in braided categories.
Lecture notes in pure and applied mathematics advances in Hopf
algebras, Vol. 158, edited by J. Bergen and S. Montgomery,   1995.

  \bibitem {Ma95b} S. Majid, Foundations of quantum group theory,
  Cambridge University Press, Cambridge, 1995.

\bibitem {Mo93}  S. Montgomery, Hopf algebras and their actions on rings.
CBMS Number 82, AMS, Providence, RI, 1993.

\bibitem{Ra85}   D.E. Radford. The structure of Hopf algebras
 with a projection.  J. Alg. {\bf 92} (1985), 322--347.

\bibitem{RT93} D.E.Radford, J.Towber,
 Yetter-Drinfeld categories associated to an arbitrary bialgebra,
 J. Pure and Applied Alg.  {\bf 87} (1993), 259-279.

\bibitem {Zh03} Shouchuan Zhang, Duality theorem  and Drinfeld double  in braided
tensor categories,  Algebra Collog., 10(2003)2, 127-134.



\end {thebibliography}
\end{document}